\documentclass[11pt]{amsart}

\allowdisplaybreaks[4]
\usepackage{graphicx} 
\usepackage{epsfig}
\usepackage{cite}
\usepackage{xypic}
\usepackage{latexsym,amssymb,amsxtra,amsmath,amsfonts}
\usepackage{verbatim}
\usepackage{color}
\usepackage{kpfonts}

\newtheorem{theorem}{Theorem}
\newtheorem{lemma}[theorem]{Lemma}
\newtheorem{proposition}[theorem]{Proposition}

\theoremstyle{definition}

\newtheorem{remark}[theorem]{Remark}

\newcommand{\beq}{\begin{equation}}
\newcommand{\eeq}{\end{equation}}
\newcommand{\beqa}{\begin{eqnarray}}
\newcommand{\eeqa}{\end{eqnarray}}
\newcommand{\beaa}{\begin{eqnarray*}}
\newcommand{\ben}{\begin{eqnarray*}}
\newcommand{\eaa}{\end{eqnarray*}}
\newcommand{\een}{\end{eqnarray*}}

\def\A{\mathcal{A}}
\def\D{\mathcal{D}}

\def\M{\mathcal{M}}
\def\H{\mathcal{H}}
\def\E{\mathcal{E}}
\def\W{\mathcal{W}}
\def\O{\mathcal{O}}
\def\F{\mathcal{F}}

\def\QQ{\mathbb{Q}}
\def\ZZ{\mathbb{Z}}

\def\CC{\mathbb{C}}
\def\PP{\mathbb{P}}
\def\ff{\mathbf{f}}
\def\ii{\mathbf{i}}
\def\tt{\mathbf{t}}
\def\qq{\mathbf{q}}

\def\lieg{{\mathfrak{g}}}
\def\normord{ {\scriptstyle {{\bullet}\atop{\bullet}}} }

\begin{document}

\title[Extended HQEs]
{\textbf{Hirota Quadratic Equations for the Gromov--Witten invariants of $\mathbb{P}_{n-2,2,2}^1$}}
\author{Jipeng Cheng}
\address{School of Mathematics, China University of Mining and
  technology, Xuzhou, Jiangsu 221116, P.R. China}
\email{chengjp@cumt.edu.cn}

\author{Todor Milanov}
\address{Kavli IPMU (WPI), UTIAS, The University of Tokyo, Kashiwa, Chiba 277-8583, Japan}
\email{todor.milanov@ipmu.jp} 

\begin{abstract}
Fano orbifold lines are classified by the Dynkin diagrams of type
$A,D,$ and $E$. It is known that the corresponding total descendant potential
is a tau-function of an appropriate
Kac--Wakimoto hierarchy. It is also known that in the A-case the Kac--Wakimoto hierarchies
admit an extension and that the total descendant potential is a
tau-function of an extended Kac--Wakimoto hierarchy. The goal of this
paper is to prove that in the D-case the total descendent potential
is also a tau-function of an extended Kac--Wakimoto hierarchy. 
\end{abstract}
\maketitle

\section{Introduction}
There are several ways to give a definition of an integrable
hierarchy. In this paper we define an integrable hierarchy in the form
of Hirota Bilinear Equations (HBEs). This approach to integrable
hierarchies is closely related to representation theory of infinite
dimensional Lie algebras and the outcome is an infinite system of
quadratic equations involving the partial derivatives of a function
$\tau$, called {\em tau-function}. The HBEs usually
have a simple geometric interpretation, i.e., they provide the
equations that define the orbit of the highest weight vector in some
highest weight representation of a given infinite dimensional Lie
algebra. On the other hand, Givental was able to reformulate the HBEs
of the $n$-KdV hierarchy in terms of the Frobenius structure on the
base space of a miniversal unfolding of the $A_{n-1}$-singularity (see
\cite{G3}). All ingredients in Givental's interpretation seem to be
generalizable to any semi-simple Frobenius manifold. A natural
question is whether HBEs can be constructed systematically for any
semi-simple Frobenius manifold and whether such a class of HBEs, if it
exists, has a representation-theoretic explanation.  So far the
answer is known to be positive for the Frobenius structures
corresponding to simple singularities  of type ADE (see \cite{GM}) and the quantum
cohomology of $\PP^1_{a,b}$ (see \cite{MT}). In general, it might be
too optimistic to expect that HBEs can be constructed for any orbifold $X$
with semi-simple quantum cohomology. Especially if the dimension of
$X$ is $>1$, then some new ideas are necessary.

On the other hand, if the target orbifold $X$ has dimension 1, then
$X$ has semi-simple quantum cohomology if and only if $X$ is an
orbifold projective line $\PP^1_{a_1,\dots,a_k}$. The 
problem of constructing HBEs seems to be approachable with the theory
of generalized Kac--Moody Lie algebras. Let us separate the orbifold
lines into three groups depending on whether the orbifold Euler
characteristic is $>0$, $=0$, or $<0$. The easiest case is the case
$>0$. The corresponding orbifold lines are called {\em Fano orbifold
  lines} and they are classified by the Dynkin diagrams of type
ADE. The problem of constructing HBEs for Fano orbifold lines was
almost completely solved in \cite{MST}. Namely, Milanov--Shen--Tseng
have constructed a system of HBEs and identified them with an
appropriate Kac--Wakimoto hierarchy, but the system of HBEs is not
complete in a sense that certain dynamical variables were
fixed. Therefore, the problem left is to construct an extension of the
Kac--Wakimoto hierarchy.

The extension in the A-case is completely understood. Both the HBEs and
the Lax equations of the corresponding integrable hierarchy are known
(see \cite{MT} for more details and other references). In this paper,
we would like to construct the extension of the Kac--Wakimoto
hierarchy of type $D$ and prove that the generating function of
Gromov--Witten (GW) invariants of $\PP^1_{n-2,2,2}$ is a solution. The precise formulation of our
result is given in Section \ref{sec:hqe}.  We also transform the HBEs
in a form convenient for the applications to integrable
hierarchies. Namely, in Section \ref{sec:change} we make a linear
change of the dynamical variables and we work out the explicit form of
our HBEs. In a subsequent paper \cite{CM}, we prove that the HBEs parametrize the
solutions to an integrable hierarchy of Lax equations, which we
propose to be called the {\em Extended D-Toda Hierarchy}.  Our two
papers solve the problem of finding an explicit description
of the integrable hierarchy that governs the GW invariants of the Fano
orbifold lines of type $D$, that is, the orbifold line
$\PP^1_{n-2,2,2}$.

\medskip

{\bf Acknowledgements.}
The work of T.M. is partially supported by JSPS Grant-In-Aid (Kiban C)
17K05193 and by the World Premier International Research Center
Initiative (WPI Initiative),  MEXT, Japan. Cheng thanks the hospitality of Kavli IPMU on his visit.

\section{Hirota quadratic equations of $\PP^1_{n-2,2,2}$ }\label{sec:hqe}

Let us recall Givental's construction of vertex operators in the settings
of quantum cohomology of the orbifold $\PP^1_{n-2,2,2}$.
\subsection{The orbifold line $\PP^1_{n-2,2,2}$}
Put  
\ben
G=\{t\in (\CC^*)^3\ |\ t_1^{n-2}=t_2^2=t_3^2\}
\een
then the orbifold $\PP^1_{n-2,2,2}$ is defined to be the translation
groupoid $[Z/G]$ where 
\ben
Z:=\{z_1^{n-2}+z_2^2+z_3^2=0\}\subset \CC^3\setminus{\{0\}}.
\een
The orbit space of $\PP^1_{n-2,2,2}$ is $Z/G\cong \PP^1$ where the
isomorphism is induced from the map
\ben
Z\to \PP^1,\quad (z_1,z_2,z_3)\mapsto [z_1^{n-2}:z_2^2].
\een
By definition the inertia orbifold
\ben
I\PP^1_{n-2,2,2}=\bigsqcup_{g\in G} [Z^g/G],\quad Z^g:=\{z\in Z\ |\ g z=z\}.
\een
The connected components of the inertia orbifold are non-empty only if
\ben
g\in G_1\cup G_2\cup G_3,
\een
where $G_i$ is the subgroup of $G$ consisting of $t=(t_1,t_2,t_3)$
such that $t_j=1$ for all $j\neq i$. Note that $G_i$ is a cyclic group
of order $a_i$, where $a_1=n-2$ and  $a_2=a_3=2.$
The connected component for $g=1$ is $\PP^1_{n-2,2,2}$, while for
$g\neq 1$ it is isomorphic to the orbifold point $[{\rm pt}/G_i]$ where $i$ is such
that $g\in G_i$. The connected components for $g\neq 1$ are
known as {\em twisted sectors}.

\subsection{Gromov--Witten invariants}
Let $H:=H^*_{\rm CR}(\PP^1_{n-2,2,2};\CC)$ be the Chen--Ruan
cohomology, i.e., 
\ben
H=\bigoplus_{g\in G_1\cup G_2\cup G_3} H^{2(*-\iota(g))}(Z^g/G,\CC), 
\een
where $*$ denotes complex degree (i.e. half of the standard
degree) and the shift $\iota(g)$ is defined for all finite order elements
\ben
g=(e^{2\pi\,  \ii \,\alpha_1}, e^{2\pi\,  \ii \,\alpha_2},
e^{2\pi\,  \ii\, \alpha_3})\in G,
\quad
\alpha_i\in \QQ\cap [0,1)
\een
by
\ben
\iota(g):=\alpha_1+\alpha_2+\alpha_3.
\een
In other words, if $\phi\in H^{2p}(Z^g/G,\CC)$, then its Chen--Ruan
degree is by definition $\operatorname{deg}_{\rm CR}(\phi)= p +\iota(g)$.
Following \cite{MST} we fix a basis of $H$ as follows
\ben
\phi_{0,0}=1,\quad \phi_{0,1}=P,
\een
and
\ben
\phi_{i,p}=1\in H^0(Z^{g_{i,p}}/G,\CC),\quad 1\leq i\leq 3,\quad 1\leq
p\leq a_i-1,
\een
where $P\in H^2(\PP^1,\CC)$ is the hyperplane class and $g_{i,p}\in
G_i$ is the element whose $i$th entry is 
$e^{2\pi\ii\, p/a_i}$. Let us define $a_0=1$, then the Chen--Ruan
degree of the above basis is given by
\ben
\operatorname{deg}_{\rm CR}(\phi_{i,p})=p/a_i,\quad 0\leq
i\leq 3.
\een
Recall that the descendant orbifold Gromov--Witten invariants of
$\PP^1_{n-2,2,2}$ are defined as intersection numbers
\ben
\langle \phi_1\psi_1^{k_1},\dots,\phi_r\psi_r^{k_r}\rangle_{g,r,d} :=
\int_{[\overline{\M}_{g,r}(\PP^1_{n-2,2,2},d)]^{\rm virt}}
\operatorname{ev}^*(\phi_1\otimes \cdots\otimes \phi_r) \psi_1^{k_1}\cdots \psi_r^{k_r},
\een
where $\phi_i\in H$, $g,r,d\in \ZZ_{\geq 0}$, and the integral is interpreted as cap product with
the virtual fundamental cycle on the moduli space
$\overline{\M}_{g,r}(\PP^1_{n-2,2,2},d)$ of orbifold stable 
maps. For more details we refer to \cite{MST} and the references
therein. The following generating function
\ben
\D(\hbar,\tt) = \exp\Big(\sum_{g,r,d=0}^\infty
\hbar^{g-1}\frac{Q^d}{r!}
\langle \tt(\psi_1),\dots,\tt(\psi_r)\rangle_{g,n,d}\Big)
\een
is called the {\em total descendant potential}. Here $Q\in \CC^*$ is a
complex parameter called the {\em Novikov variable}, $\hbar$,
$t_0,t_1,\dots\in H$ are formal vector variables, and
$\tt(z)=\sum_{k=0}^\infty t_k z^k$. The components of the vector
variable $t_k$ with respect to the basis $\phi_{i,p}$ of $H$ from
above will be denoted by $t_{i,p,k}$.  

\subsection{The calibration operator}
If $t\in H$ then we define
\ben
\langle \phi_1\psi_1^{k_1},\dots,\phi_r\psi_r^{k_r}\rangle_{g,r}(t) =
\sum_{\ell,d=0}^\infty
\frac{Q^d}{\ell!} \langle \phi_1\psi_1^{k_1},\dots,\phi_r\psi_r^{k_r},t,\dots,t\rangle_{g,r+\ell,d}.
\een
Using the divisor equation and degree reasons one can prove that the
above correlator is polynomial in $t$ and $Qe^{t_{0,1}}$. In
particular, it defines an analytic function on $H$.

The quantum cup product  $\bullet_{t,Q}$ ($(t,Q)\in H\times \CC^*$) is a family of
multiplications in $H$ defined by the identity
\ben
(\phi_1\bullet_{t,Q} \phi_2,\phi_3):= \langle\phi_1,\phi_2,\phi_3\rangle_{0,3}(t)
\een
for all $\phi_1,\phi_2,\phi_3\in H$, where
\ben
(\phi_1,\phi_2) := \langle \phi_1,\phi_2,1\rangle_{0,3,0} 
\een
is the orbifold Poincare pairing. The specialization of the quantum
cup product to $t=Q=0$ is known as the Chen--Ruan cup product $\cup_{\rm
  CR}:=\bullet_{0,0}$.  Let us point out also that the orbifold
Poincare pairing is given explicitly by the following formulas
\ben
(\phi_{i,p},\phi_{j,q}) = \frac{1}{a_i}
\,\delta_{i,j}\delta_{p+q,a_i},\quad 0\leq i,j\leq 3.
\een
The multiplication and the Poincare pairing define a conformal Frobenius
structure on $H$ of conformal dimension 1 with Euler vector field
\ben
E=\sum_{i,p} (1-p/a_i)t_{i,p}\frac{\partial}{\partial t_{i,p}} +
\frac{1}{n-2} \frac{\partial}{\partial t_{0,1}}.
\een
Note that $\chi_{\rm orb}(\PP^1_{n-2,2,2}) = 1/(n-2)$ is the orbifold
Euler characteristic. In particular, the connection $\nabla$ on the vector
bundle $TH\times \CC^*\to H\times \CC^*$ defined by
\ben
\nabla_{\partial/\partial t_{i,p}} & = &   \frac{\partial}{\partial
  t_{i,p}} -z^{-1} \phi_{i,p}\bullet_{t,Q}\\
\nabla_{\partial/\partial z} & = &  \frac{\partial}{\partial z}
+z^{-2} E\bullet_{t,Q} - z^{-1}\theta
\een
is flat, where
\ben
\theta:H\to H,\quad \theta(\phi_{i,p}) = \Big( \frac{1}{2} -
\frac{p}{a_i}\Big) \phi_{i,p}
\een
is the so called {\em Hodge grading operator}. 
We refer to \cite{Du} for further details on Frobenius structures.

The Dubrovin's connection has a solution of the form
$S(t,Q,z) z^\theta z^{-\rho}$, where
\ben
\rho = \frac{1}{n-2} P \cup_{CR} = c_1(T\PP^1_{n-2,2,2})\cup_{\rm CR}.  
\een
is a nilpotent operator and
\ben
S(t,Q,z)=1+S_1(t,Q) z^{-1} +S_2(t,Q)z^{-2}+\cdots,\quad
S_k(t,Q)\in \operatorname{End}(H) 
\een
is an operator series defined by
\ben
(S(t,Q,z)a,b) = (a,b)+\sum_{k=0}^\infty
z^{-k-1}\langle a\psi_1^k,b\rangle_{0,2}(t).
\een
\subsection{Periods}
Let us recall the second structure connection
\ben
\nabla^{(\mu)}_{\partial_{i,p}} & = &  
\partial_{i,p} + (\lambda-E\bullet)^{-1} (\phi_{i,p} \bullet) (\theta-\mu-1/2) \\
\nabla^{(\mu)}_{\partial_\lambda} & = & 
\partial_\lambda -(\lambda-E\bullet_t)^{-1} (\theta-\mu-1/2),
\een
where $\mu\in \CC$ is a complex parameter and
$\partial_{i,p}:=\partial/\partial t_{i,p}$. 
This is a connection on the trivial bundle
\ben
(H\times \CC)'\times H \to (H\times \CC)',
\een
where 
\ben
(H\times \CC)'=\{ (t,\lambda)\ |\ \det (\lambda-E\bullet)\neq 0\}.
\een
The hypersurface $\det (\lambda-E\bullet_t)=0$ in $H\times \CC$ is
called the {\em discriminant}.

Let us fix a reference point $(t^\circ,\lambda^\circ)\in (H\times
\CC)'$, such that $\lambda^\circ$ is a sufficiently
large real number. Suppose that $m\in \ZZ$ is an integer, such that
$m<\alpha$ for all eigenvalues $\alpha$ of $\theta+\tfrac{1}{2}$. It is
easy to check that the following functions provide a fundamental
solution to the 2nd structure connection $\nabla^{(m)}$ (see \cite{Mi})
\ben
I^{(m)}(t,\lambda) = \sum_{k=0}^\infty (-1)^k S_k(t,Q) \widetilde{I}^{(m+k)}(\lambda),
\een 
where
\ben
\widetilde{I}^{(m)}(\lambda) = e^{-\rho \partial_\lambda \partial_m}
\Big(
\frac{\lambda^{\theta-m-\frac{1}{2}} }{ \Gamma(\theta-m+\frac{1}{2}) }
\Big).
\een
Note that the following relation is satisfied:
\beq\label{transl-invar}
\partial_\lambda I^{(m)}(t,\lambda) = I^{(m+1)}(t,\lambda)
\eeq
for all $m$, such that $m+1<\alpha$  for all eigenvalues $\alpha$ of
$\theta+\tfrac{1}{2}$. Therefore, we can uniquely extend the
definition of $I^{(m)}(t,\lambda)$ for all $m\in \ZZ$ in such a way
that \eqref{transl-invar} is satisfied for all $m\in \ZZ$.
The 2nd structure connection has a
Fuchsian singularity at infinity, therefore the series $I^{(m)}(t,\lambda)$ is
convergent for all $(t,\lambda)$ sufficiently close to
$(t^\circ,\lambda^\circ)$. Using the differential equations we extend
$I^{(m)}$ to a multi-valued analytic function on $(H\times \CC)'$. We
define the following multi-valued functions taking values in $H$:
\ben
I^{(m)}_a(t,\lambda):=I^{(m)}(t,\lambda)\, a, 
\quad a\in H,\quad m\in \ZZ.
\een 
The functions $I^{(m)}_a$ and $\widetilde{I}^{(m)}_a$ will be called
respectively {\em periods} and {\em calibrated periods}.

Using analytic continuation we get a representation 
\beq\label{mon-repr}
\pi_1((H\times\CC)',(t^\circ,\lambda^\circ) )\to \operatorname{GL}(H)
\eeq
called the {\em monodromy representation} of the Frobenius
manifold. The image $W$ of the monodromy representation is called the
{\em monodromy group}. 

Using the differential equations of the 2nd structure connection it is
easy to prove that the pairing
\ben
(a|b):=(I^{(0)}_a(t,\lambda),(\lambda-E\bullet)I^{(0)}_b(t,\lambda))
\een
is independent of $t$ and $\lambda$. This pairing is known as the {\em
  intersection pairing}.
The monodromy group $W$ is generated by a set of reflections 
\ben
w_a(x)=x-(a|x) a,\quad a\in \mathfrak{R},
\een
where the set $\mathfrak{R}$ is defined as follows. Since the quantum cohomology of $\mathbb{P}_{n-2,2,2}^1$
is semisimple, we may choose a generic reference point, such that 
the Frobenius multiplication $\bullet_{t^\circ,Q}$ is semi-simple and the
operator $E\bullet_{t^\circ,Q}$ has $N:=n+1$ pairwise different eigenvalues
$u_i^\circ$. Let $\mathfrak{R}$ be the set of all $a\in H$ such that $(a|a)=2$
and there exists
a simple loop in $\CC-\{u_1^\circ,\cdots,u_N^\circ\}$ based at $\lambda^\circ$
such that monodromy transformation along it transforms $a$ into
$-a$. Here simple loop means a loop that starts at $\lambda^\circ$,
approaches one of the punctures $u_i^\circ$ along a path $\gamma$
that ends at a point sufficiently close to $u_i^\circ$, goes around
$u_i^\circ$, and finally returns back to $\lambda^\circ$ along $\gamma$.

\subsection{The integral lattice of Iritani}

Let $X$ be a complex orbifold groupoid whose orbit space $|X|$ is a projective
variety. Using the $K$-ring $K^0(X)$ of topological orbifold vector
bundles on $X$ and a certain $\Gamma$-modification of the Chern
character map, Iritani has introduced an integral lattice in the
Chen-Ruan cohomology group $H_{\rm CR}(X;\CC)$. If $X$ has semi-simple
quantum cohomology, then it is expected that $X$ has a LG mirror model
and that Iritani's lattice coincide with the image of the Milnor
lattice via an appropriate period map. Let us recall Iritani's
construction in a form suitable for our purposes.

Let $IX$ be the {\em inertia orbifold} of $X$, that is, as a groupoid the points of $IX$ are 
\ben
(IX)_0=\{(x,g)\ |\ x\in X_0,\  g\in \operatorname{Aut}(x)\}
\een
while the arrows from $(x',g')$ to $(x'',g'')$ consists of all arrows
$g\in X_1$ from $x'$ to $x''$, such that, 
$g''\circ g = g\circ g'$. It is known that $IX$ is an orbifold
consisting of several connected components $X_v$, $v\in
T:=\pi_0(|IX|)$.

Let us fix an ample basis $\{P_i\}_{i=1}^r\subset H^2(|X|,\ZZ)$ and
let $Q_i$ $(1\leq i\leq r)$ be the corresponding Novikov variables.
Following Iritani \cite{Ir}, we define a linear map  
\ben
\Psi: K^0(X)\to H^*(IX;\CC)=\oplus_{v\in T} H^*(X_v;\CC)
\een
via 
\beq\label{psi-map}
\Psi(E)=(2\pi)^{(1-\operatorname{dim}  X)/2}\ 
\Big(\widehat{\Gamma}(X) e^{-\sum_{i=1}^rP_i\log Q_i}\Big)\cup
(2\pi\sqrt{-1})^{\operatorname{deg}}
\operatorname{inv}^* \widetilde{\operatorname{ch}}(E).
\eeq
Let us recall the rest of the notation. The linear operator 
\ben
\operatorname{deg}: H^*(IX;\CC)\to H^*(IX;\CC)
\een 
is defined by $\operatorname{deg} (\phi)=i \phi$ if $\phi \in H^{2i}(IX;\CC)$. 
The involution $\operatorname{inv}:IX\to IX$ inverts all arrows while
on the points it acts as $(x,g)\mapsto (x,g^{-1}).$ If $E$ is an
orbifold vector bundle, then we have an eigenbasis decomposition 
\ben
\operatorname{pr}^*(E) = 
\oplus_{v\in T} E_v=
\oplus_{v\in T}\oplus_{0\leq f<1} E_{v,f},
\een
where $\operatorname{pr}:IX\to X$ is the forgetful map $(x,g)\mapsto
x$ and $E_{v,f}$ is the subbundle of
$E_v:=\operatorname{pr}^*(E)|_{X_v}$ whose fiber over a point
$(x,g)\in (IX)_0$ is the eigenspace of $g$ corresponding to the
eigenvalue $e^{2\pi\sqrt{-1} f}$. Let us denote by $\delta_{v,f,i}$
$(1\leq i\leq l_{v,f}:={\rm rk}(E_{v,f}))$ the Chern roots of
$E_{v,f}$, then the Chern character and the $\Gamma$-class of $E$ are
defined by 
\ben
\widetilde{\operatorname{ch}}(E) = 
\sum_{v\in T} \sum_{0\leq f<1} e^{2\pi\sqrt{-1} f}
\sum_{i=1}^{l_{v,f}} e^{\delta_{v,f,i}}
\een
\ben
\widehat{\Gamma}(E) = 
\sum_{v\in T} 
\prod_{0\leq f<1} 
\prod_{i=1}^{l_{v,f}} 
\Gamma(1-f+\delta_{v,f,i}),
\een
where the value of the $\Gamma$-function $\Gamma(1-f+y)$  at
$y=\delta_{v,f,i}$ is obtained by first expanding in Taylor's series
at $y=0$ and then formally substituting $y=\delta_{v,f,i}$. By
definition $\widehat{\Gamma}(X):=\widehat{\Gamma}(TX)$. The cup
product in \eqref{psi-map} is the usual topological cup product on
$|IX|$.   
\begin{remark}
There are two differences between the maps $\Psi$ defined respectively by
\eqref{psi-map}  and formula (37) in \cite{MST}. First, the
normalization factor in \cite{MST} is $(2\pi)^{-\operatorname{dim}  X/2}$ instead of
$(2\pi)^{(1-\operatorname{dim}  X)/2}$ and second, the map $\Psi$ in
\cite{MST} does not depend on the Novikov variables.\qed
\end{remark}
\begin{lemma}\label{le:ip}
  The intersection pairing is independent of the Novikov variable $Q$ and the
  following formula holds:
  \ben
  (a|b) = \langle a,b\rangle+\langle b,a\rangle, \quad
  \forall a,b\in H,
  \een
  where
  \ben
  \langle a,b\rangle := \frac{1}{2\pi} \left(a,
    e^{\pi\sqrt{-1}\theta}
    e^{\pi\sqrt{-1}\rho} b\right).
  \een
\end{lemma}
The proof of Lemma \ref{le:ip} will be given in Section \ref{sec:ip}. 
Using the Kawasaki--Riemann--Roch formula, one can 
prove that the map $\Psi$ intertwines the pairing
$\langle\ ,\ \rangle$ defined in 
Lemma \ref{le:ip} and the Euler pairing on $K^0(X)$, that is,
\beq\label{euler-pairing}
\chi(E_1^\vee\otimes E_2) =
\langle \Psi(E_1),  \Psi(E_2)\rangle .
\eeq
We refer to \cite{Ir}, Proposition 2.10, (iii) for the details of this
computation.
\begin{remark}
Formula \eqref{euler-pairing} is used only to
conclude that its RHS is an integer, which in our settings could also
be checked directly, because we will compute the RHS of
\eqref{euler-pairing} explicitly.
\end{remark}
According to Lemma \ref{le:ip}, we have
\beq\label{int=sym-euler}
(\Psi(a)|\Psi(b)) = \chi(a^\vee\otimes b)+\chi(a\otimes b^\vee),\quad
a,b\in K^0(X).
\eeq
From now on we will use $\Psi$ to identify $H$ and $K^0(X)\otimes
\CC$. In particular, the above formula allows us to say that the
intersection pairing is the symmetrization of the Euler pairing.

Let us specialize again to $X=\PP^1_{n-2,2,2}$. In this case, we choose $P_1$ to be
the hyperplane class $P$. 
According to Milanov--Shen--Tseng (see Theorem 12 and
Proposition 13 in \cite{MST}) the set of reflection vectors
\ben
\mathfrak{R}= \{\Psi(a)\ |\ a\in K^0(X) \mbox{ such that } (a|a)=2\}.
\een

\subsection{Calibrated periods}
The $K$-ring of $\PP^1_{n-2,2,2}$ can be identified with
\ben
K:= \ZZ[L_1,L_2,L_3]/\langle L_i^{a_i}-L_j^{a_j}, (L_i-1)(L_j-1)\ |\
1\leq i\neq j\leq 3\rangle.
\een
The generators $L_i$ correspond to the orbifold line bundles
$[(Z\times \CC)/G]$, where $G$ acts on the fiber $\CC$ via the
character
\ben
G\to \CC^*, \quad g=(g_1,g_2,g_3)\mapsto g_i.  
\een
The orbifold tangent bundle $T\PP^1_{n-2,2,2}$ can be identified with $[Q/G]$
where $Q$ is the $G$-equivariant bundle on $Z$ defined as a quotient
\ben
0\to Z\times \lieg \to TZ \to Q\to 0,
\een
where $\lieg$ is the Lie algebra of $G$ and the map $Z\times \lieg \to
TZ$ is defined by
\ben
(z,\xi)\mapsto \left.\frac{d}{d\epsilon} \Big(
e^{\epsilon \xi} z\Big)\right|_{\epsilon=0} \in T_zZ. 
\een
Note that in the $K$-ring we have
\ben
[TZ/G]=L_1+L_2+L_3-L,
\een
where $L:=L_i^{a_i}$. Therefore
\ben
T\PP^1_{n-2,2,2} = L_1+L_2+L_3-L-1 = L_1L_2L_3L^{-1}. 
\een

Let us introduce the following notation for the components of the
Chern character map:
\ben
\widetilde{\operatorname{ch}}(E) = 
\operatorname{rk}(E)\phi_{0,0}+ 
\operatorname{deg}(E)\phi_{0,1} + \sum_{j=1}^3\sum_{p=1}^{a_j-1}
\chi_{j,p}(E) \phi_{j,p}.
\een
The first component $\operatorname{rk}(E)$ is the rank of 
$E$, the second one $\operatorname{deg}(E) :=\int_{|\PP^1_{n-2,2,2}|} c_1(E)$ is the
degree of $E$, and the remaining ones 
$\chi_{j,p}:K\to \CC$ are some {\em ring} homomorphisms. The rank is
straightforward to compute, while for the remaining components we have
the following explicit formulas:
\ben
\operatorname{deg}(L_i)=1/a_i,\quad 
\chi_{j,p}(L_i) :=e^{-\frac{2\pi\ii}{a_j} p \delta_{j,i}}.
\een
Recalling the definition, we get that Iritani's map 
$\Psi$ has the following explicit form (see also \cite{MST},
formula (39)):
\ben
\Psi(E) = \operatorname{rk}(E)\Big(1 -\frac{\gamma}{n-2}P\Big) +2\pi\ii
  \operatorname{deg}(E) P +
\sum_{j=1}^3\sum_{p=1}^{a_j-1} \Gamma(1-p/a_j) \chi_{j,p}(E) \phi_{j,p},
\een
where $\gamma=-\Gamma'(1)+(n-2)\log Q$.  

Using formulas \eqref{euler-pairing} and \eqref{int=sym-euler},
we get (see Section \ref{sec:eup} for similar computations) 
\ben
(L_i^m|L_j^k)=0\quad
\mbox{for all } \quad
i\neq j ,\quad
1\leq m\leq a_i-1,\quad 
1\leq k \leq a_j-1
\een
and
\ben
(L_i^m|L_i^k)=
\begin{cases}
  2 & \mbox{if } m\equiv k({\rm mod }\  a_i),\\
  1 & \mbox{otherwise}.
\end{cases}
\een
The pairing $(\ |\ )$ is degenerate with kernel spanned by
$(L-1)$. Using the above formulas or the formulas in Section
\ref{sec:eup}, we get that the following vectors
project to an orthonormal basis of $K/\ZZ(L-1)$:
\ben
\epsilon^1_i & = & L_1^i+\frac{1}{2} (L_2+L_3)-1\quad
(1\leq i\leq n-2) ,\\
\epsilon^2_1 & = & \frac{1}{2} (L_2+L_3)-1, \\
\epsilon^3_1 & = & \frac{1}{2} (L_2-L_3).
\een
\begin{proposition}\label{pr:cal-periods}
a) The calibrated period corresponding to an arbitrary $\alpha\in K$ is
given by  the following formula
\ben
\widetilde{I}^{(-\ell-1)}_\alpha(\lambda)  &=
&
\operatorname{rk}(\alpha) \frac{\lambda^{\ell+1}}{(\ell+1)!} \, \phi_{0,0}+
\frac{\lambda^{\ell}}{\ell!} \Big(
\tfrac{1 }{n-2} \, \operatorname{rk}(\alpha) (\log
\lambda -C_\ell) +2\pi\ii \operatorname{deg}(\alpha)\Big)\, \phi_{0,1} + \\
&&
\sum_{j=1}^3\sum_{p=1}^{a_j-1}
\frac{\lambda^{\ell+1-p/a_j}}
{\Big(\ell+1-\frac{p}{a_j}\Big)\cdots
  \Big(1-\frac{p}{a_j}\Big) }\,
\chi_{j,p}(\alpha)\, \phi_{j,p},
\een
where $C_0:=(n-2)\log Q$ and  $C_{\ell+1}=C_{\ell}+\frac{1}{\ell+1}$ for
$\ell\geq 0$.

b) Let $\sigma:K\to K$ be the classical monodromy operator defined by:
the analytic continuation in
anti-clockwise direction of $ \widetilde{I}^{(m)}_\alpha(\lambda)$ is
$ \widetilde{I}^{(m)}_{\sigma(\alpha)}(\lambda)$. Then $\sigma(\alpha)
= \alpha\, T\PP^1_{n-2,2,2}$.

c) We have
\ben
\sigma(\epsilon^1_i)  & = & 
\begin{cases}
  \epsilon^1_{i+1}, & 1\leq i\leq n-3 \\
  \epsilon^1_1+L-1, & i=n-2,
\end{cases}\\
\sigma(\epsilon^2_1)  & = & -\epsilon^2_1+L-1, \\
\sigma(\epsilon^3_1)  & = & -\epsilon^3_1.
\een
\end{proposition}
The proof of Proposition \ref{pr:cal-periods} is a straightforward
computation using formula \eqref{calibr-period}, so it will be omitted. 
Finally,  let us point out that the set of reflection vectors is given
explicitly by
\beq\label{refl_vect}
\mathfrak{R} = \{\pm(\epsilon^a_i\pm\epsilon^b_j)+m(L-1)\},
\eeq
where $(a,i)\neq (b,j)$ and $m\in \ZZ$ take all possible values.

\subsection{Vertex operators}

Let us change the variables in the total descendant potential
\ben
q_{k}:=t_{k}-\delta_{k,1}\mathbf{1},\quad 
k\in \ZZ_{\geq 0},
\een
where $\mathbf{1}:=\phi_{0,0}$ is the unit. The components of the formal
vector variable $q_k$ with respect to the basis $\{\phi_{i,p}\}\subset
H$ will be denoted by $q_{i,p,k}$. This substitution identifies the
total descendant potential with a vector in the ring of formal power
series with a shifted origin
$\CC_\hbar[\![q_0,q_1+\mathbf{1},q_2,\dots]\!]$, where
$\CC_\hbar=\CC(\!(\hbar^{1/2})\!)$.   

Following Givental \cite{G2} we introduce the symplectic vector space
$\H:=H(\!(z^{-1})\!)$ with symplectic form
\ben
\Omega(f,g)=\operatorname{Res}_{z=0} (f(-z),g(z)) dz.
\een
The vector space $\H\oplus \CC$ has a natural structure of a Heisenberg
Lie algebra
\ben
[a,b]:=\Omega(a,b),\quad a,b\in \H.
\een
Let $\{\phi^{i,p}\}\subset H$ be the basis dual to $\{\phi_{i,p}\}$
with respect to the orbifold Poincare pairing.  
The formulas
\ben
(\phi_{i,p} z^k)\sphat:= -\sqrt{\hbar}\partial/\partial q_{i,p,k},\quad 
(\phi^{i,p}(-z)^{-k-1})\sphat:=q_{i,p,k}/\sqrt{\hbar},
\een
define a representation of the Heisenberg Lie algebra on
$\CC_\hbar[\![q_0,q_1+\mathbf{1},q_2,\dots]\!]$.  

Given $\alpha\in K$ put 
\ben
\ff_\alpha(t,\lambda,z)=\sum_{m\in \ZZ} I^{(m)}_\alpha(t,\lambda) (-z)^m.
\een
and 
\ben
\widetilde{\ff}_\alpha(\lambda,z) = 
\sum_{m\in \ZZ} 
\widetilde{I}^{(m)}_\alpha(\lambda) (-z)^m.
\een
We will be interested in the vertex operators 
\ben
\Gamma^\alpha(t,\lambda):= 
e^{\ff^-_\alpha(t,\lambda,z)\sphat}
e^{\ff^+_\alpha(t,\lambda,z)\sphat}
\een
and 
\ben
\widetilde{\Gamma}^\alpha(\lambda):= 
e^{\widetilde{\ff}^-_\alpha(\lambda,z)\sphat}
e^{\widetilde{\ff}^+_\alpha(\lambda,z)\sphat},
\een
where the superscript $+$ (resp. $-$) denotes the $z$-series obtained by
truncating all terms that contain negative (resp. non-negative) powers
of $z$. 

Let us denote by $\A$ the algebra of differential operators in one
variable
\ben
A(x,\partial_x)=\sum_{i=0}^d A_i(x)\partial_x^i,\quad 
A_i(x)\in \CC[x].
\een
The algebra is equipped with an anti-involution 
\ben
{}^\#:\A\to \A,\quad 
A\mapsto A^\#:=\sum_{i=0}^d (-\partial_x)^i \circ A_i(x),
\een
where $\circ$ means composition (not action!).  
Let us introduce the following vertex operator with coefficients in
$\A$:
\ben
\widetilde{\Gamma}(\lambda) = 
\exp\Big(
-\sum_{\ell>0} \frac{\lambda^\ell}{\ell!} \partial_x q_{0,0,\ell}\Big)
\exp\Big(x\partial_{0,0,0}\Big),
\een
where $\partial_{i,p,k}:=\partial/\partial q_{i,p,k}$. 

Now we can introduce the main object of our investigation. That is a
system of quadratic equations, which will be called the Hirota
Quadratic Equations (HQEs) of $\PP^1_{n-2,2,2}$.  
Put
\ben
\E:=\{\pm \epsilon^1_i\ (1\leq i\leq n-2),\ \pm \epsilon^2_1, \ \pm \epsilon^3_1\},
\een
\ben
\widetilde{b}_{\pm \epsilon^1_i}(\lambda) = \frac{Q}{n-2} \lambda^{-1/(n-2)}
\eta_1^{-i},
\quad
\eta_1:=e^{2\pi\ii /(n-2)},
\een
and 
\ben
\widetilde{b}_{\pm \epsilon^2_1}(\lambda) =-\frac{1}{4},\quad 
\widetilde{b}_{\pm \epsilon^3_1}(\lambda)=\frac{1}{4}.   
\een
The coefficients $\widetilde{b}_\epsilon$ can be defined also in terms
of the phase factors corresponding to the composition of the vertex
operators $\widetilde{\Gamma}^\epsilon(\lambda)
\widetilde{\Gamma}^{-\epsilon}(\mu)$ (see Section \ref{sec:pf} and Lemma
\ref{le:calibr-b} for more details).

We say that a function (or formal power series in $\mathbf{q}=(q_{i,p,k})$) $\tau(\hbar,\qq)$
satisfies the HQEs of $\PP^1_{n-2,2,2}$ if for
every integer $m\in \ZZ$ the 1-form 
\beq\label{hqe}
\frac{d\lambda}{\lambda}
\Big(\widetilde{\Gamma}^\#(\lambda)\otimes \widetilde{\Gamma}(\lambda)\Big)
\Big(\sum_{\epsilon\in \E} 
\widetilde{b}_{\epsilon}(\lambda) 
\widetilde{\Gamma}^\epsilon(\lambda)\otimes
\widetilde{\Gamma}^{-\epsilon}(\lambda)
\Big)(\tau\otimes \tau)
\eeq
computed at $q_{0,0,0}'-q_{0,0,0}''=m\hbar^{1/2} $ is regular in
$\lambda$. Note that the entire expression makes sense as an element
in 
\ben
d\lambda 
\A(\!(\lambda^{-1})\!)(\!(\hbar^{1/2})\!)[\![ \qq'+\mathbf{1},\qq''+\mathbf{1}]\!].
\een
The  regularity requirement means that the coefficients in front of
the monomials involving $\qq'$, $\qq''$, and $\hbar^{1/2}$ are in fact
polynomial in $\lambda$.  
\begin{theorem}\label{t1}
The total descendant potential of $\PP^1_{n-2,2,2}$ satisfies the
HQEs of $\PP^1_{n-2,2,2}$. 
\end{theorem}

Let us clarify the relation of our HQEs
\eqref{hqe} to the Kac--Wakimoto hierarchies. To begin with, note that
the set of reflection vectors \eqref{refl_vect}, modulo the rank 1
lattice $\ZZ (L-1)$, 
form a root system of type $D$ with respect to the intersection
pairing. It is an easy exercise, using the formulas in Proposition
\ref{pr:cal-periods} c), to check that the automorphism $\sigma$
induces an  
element of the Weyl group given by the composition of the
reflections corresponding to the non-branching nodes of the Dynkin
diagram (see also \cite{MST}, Proposition 17).  If we specialize
$q_{0,0,\ell}\otimes 1=1\otimes q_{0,0,\ell}$ for all $\ell>0$, then
\eqref{hqe} becomes a system 
of HBEs which is equivalent to the HBEs of the Kac--Wakimoto 
hierarchy of type $D$ corresponding to the conjugacy class of
$\sigma$. Let us point out that here we get a realization of the
Kac--Wakimoto hierarchy corresponding to a fermionic realization of
the basic representation (see \cite{KvL}). 

\begin{remark}
  Our interpretation of \eqref{hqe} leads to a system of
  quadratic equations for the Taylor's coefficients of the
  tau-function $\tau$. On the other hand, in the theory of integrable
  systems, there is a different interpretation that leads to a system
  of PDEs, which is usually called HBEs (e.g. see \cite{Kac}, Section 14.11).
  We refer to \eqref{hqe} as HQEs
  or HBEs depending on whether we would like to think of \eqref{hqe}
  as a system of algebraic equations or a system of PDEs.
\end{remark}

\section{The extended D-Toda hierarchy}\label{sec:change}

In this section, by making an explicit linear change of the variables
$q_{i,p,k}$, we will transform the HQEs \eqref{hqe} into form
convenient for the applications to integrable 
hierarchies. In a companion paper \cite{CM} to this one, we prove that our
system of HBEs  parametrizes the solutions to an
integrable hierarchy, which we suggest to be called the {\em Extended
  D-Toda Hierarchy}. 
\subsection{Change of variables}
Let $\tt=(\tt_0,\tt_1,\tt_2,\tt_3)$ be 4 sequences of formal variables,
where
\ben
\tt_0=(t_{0,\ell})_{\ell\geq 1},\quad
\tt_1=(t_{1,\ell})_{\ell\geq 1},\quad
\tt_2=(t_{2,2\ell+1})_{\ell\geq 0},\quad
\tt_3=(t_{3,2\ell+1})_{\ell\geq 0}.
\een
Let us introduce the vertex operators
\ben
\Gamma_1^\pm(z) & := & \exp\Big(\pm \sum_{\ell=1}^\infty t_{1,\ell} z^\ell\Big)
\exp\Big(\mp\sum_{\ell=1}^\infty \frac{z^{-\ell}}{\ell}
                  \frac{\partial}{\partial t_{1,\ell}}\Big),\\
\Gamma_a(z) & := &\exp\Big( \sum_{\ell=0}^\infty t_{2,2\ell+1} z^{2\ell+1}\Big)
\exp\Big(-2\sum_{\ell=0}^\infty \frac{z^{-2\ell-1}}{2\ell+1}
\frac{\partial}{\partial t_{a,2\ell+1}}\Big),\quad
a=2,3.
\een
Using Proposition \ref{pr:cal-periods} we get the following formulas
\ben
&&
\widetilde{\Gamma}^{\epsilon^1_i}(\lambda)=
\exp\left(
\sum_{\ell=0}^\infty \frac{\lambda^\ell}{\ell!}
\Big(\frac{1}{n-2}(\log\lambda-C_\ell)+2\pi\ii\Big(\frac{i}{n-2}+\frac{1}{2}\Big)\Big)
q_{0,0,\ell}/\sqrt{\hbar}\right)
e^{-\sqrt{\hbar}\partial_{0,0,0}}\\
&&
\exp\left(
\sum_{\ell=0}^\infty\Big( \frac{\lambda^{\ell+1}}{(\ell+1)!}\,
q_{0,1,\ell}/\sqrt{\hbar}+
\sum_{j=1}^{n-3}
\frac{\eta_1^{ij}\lambda^{\ell+\frac{j}{n-2}}  }{
\Big(\ell+\frac{j}{n-2}\Big)\cdots\Big(0+\frac{j}{n-2}\Big)(n-2)}
q_{1,j,\ell}/\sqrt{\hbar}
\Big)\right)\\
&&
\exp\left(
  -\sum_{\ell=0}^\infty\Big(
  \frac{\ell!}{n-2}\lambda^{-1-\ell}\sqrt{\hbar}\partial_{0,1,\ell}+
  \sum_{j=1}^{n-3}\prod_{k=0}^{\ell-1}
  \Big(k+\frac{j}{n-2}\Big)
  \lambda^{-\frac{j}{n-2}-\ell}\eta_1^{-ij}\sqrt{\hbar}\partial_{1,j,\ell}
  \Big)\right).
\een
Let us make the following substitutions: $\lambda
=\frac{z_1^{n-2}}{n-2}$, 
\beq\label{subst-1k1}
t_{1,\ell(n-2)}
=\hbar^{-1/2}
\frac{q_{0,1,\ell-1}}{(n-2)^{\ell}\ell!},\quad \ell> 0, 
\eeq
\beq\label{subst-1k2}
  t_{1,\ell(n-2)+i} =\hbar^{-1/2}\, (n-2)^{-\frac{i}{n-2}}
  \frac{
  q_{1,i,\ell}}{
  i(i+n-2)\cdots(i+(n-2)\ell) },\quad 1\leq i\leq n-3,\quad \ell \geq 0.
\eeq
Then the vertex operator
\ben
\widetilde{\Gamma}^{\pm\epsilon^1_i}(\lambda) =
e^{\pm
\sum_{\ell=0}^\infty \frac{\lambda^\ell}{\ell!}
\Big(\frac{1}{n-2}(\log\lambda-C_\ell)+
2\pi\ii\left(\frac{i}{n-2}+\frac{1}{2}\right)\Big)
q_{0,0,\ell}/\sqrt{\hbar}}
e^{\mp\sqrt{\hbar}\partial_{0,0,0}}
\Gamma_1^{\pm}(\eta_1^i z_1).
\een
and the 1-form
\ben
\widetilde{b}_{\epsilon^1_i}(\lambda)\frac{d\lambda}{\lambda}=
-C\eta_1^{-i} \frac{dz_1}{z_1^2},\quad
C:=-Q(n-2)^{\frac{1}{n-2}}.
\een
A straightforward computation yields
\ben
\widetilde{\Gamma}^\#\widetilde{\Gamma}^{\pm \epsilon^1_i}\tau =
\Big(\Gamma^\pm_1(\eta_1^i z_1)\tau(x\mp\sqrt\hbar,\qq)\Big)
e^{
  \sum_{\ell>0} \frac{\lambda^\ell}{\ell!}
  (\sqrt\hbar \partial_x \mp h_\ell)q_{0,0,\ell}/\sqrt{\hbar}}
\Big(\frac{\eta_1^i z_1}{C}\Big)^
{\pm(q_{0,0,0}+x)/\sqrt\hbar}
\een
and
\ben
\widetilde{\Gamma}\widetilde{\Gamma}^{\pm \epsilon^1_i}\tau =
\Big(\frac{\eta_1^i z_1}{C}\Big)^
{\pm(q_{0,0,0}+x)/\sqrt\hbar}
e^{
  \sum_{\ell>0} \frac{\lambda^\ell}{\ell!}
  (-\sqrt\hbar \partial_x \mp h_\ell)q_{0,0,\ell}/\sqrt{\hbar}}
\Big(\Gamma^\pm_1(\eta_1^i z_1)\tau(x\mp\sqrt\hbar,\qq)\Big),
\een
where $\tau(x,\qq):=e^{x\partial_{0,0,0}}\tau(\qq)$ and
$h_\ell=\tfrac{1}{n-2}\left(1+\tfrac{1}{2}+\cdots+\tfrac{1}{\ell}\right)$
for $\ell>0$.  The residue in the $\lambda$-plane
\ben
\operatorname{Res}_{\lambda=\infty}
\frac{\lambda^r}{r!}\,
\sum_{i=1}^{n-2} \widetilde{b}_{\epsilon^1_i}(\lambda)
\frac{d\lambda}{\lambda}
\Big(
\widetilde{\Gamma}^\#\otimes \widetilde{\Gamma}
\Big)\Big(
\widetilde{\Gamma}^{\epsilon^1_i}(\lambda)\otimes
\widetilde{\Gamma}^{-\epsilon^1_i}(\lambda) +
\widetilde{\Gamma}^{-\epsilon^1_i}(\lambda)\otimes
\widetilde{\Gamma}^{\epsilon^1_i}(\lambda)
\Big)(\tau\otimes\tau)
\een
turns into the following residue in the $z_1$-plane
\ben
&&
-\operatorname{Res}_{z_1=\infty}
\frac{z_1^{(n-2)r}}{(n-2)^rr!}
\frac{dz_1}{z_1}\left(
  (z_1/C)^{m-1}\times
  \phantom{
    e^{
  \sum_{\ell>0} \frac{z_1^{(n-2)\ell}}{(n-2)^\ell\ell!}
  (\sqrt\hbar \partial_x -
  h_\ell)(q'_{0,0,\ell}-q''_{0,0,\ell})/\sqrt{\hbar}}
}
\right.\\
&&
\Big(\Gamma^+_1(z_1)\tau(x-\sqrt\hbar,\qq')\Big)
e^{
  \sum_{\ell>0} \frac{z_1^{(n-2)\ell}}{(n-2)^\ell\ell!}
  (\sqrt\hbar \partial_x -
  h_\ell)(q'_{0,0,\ell}-q''_{0,0,\ell})/\sqrt{\hbar}}
\Big(\Gamma^-_1(z_1)\tau(x+\sqrt\hbar,\qq'')\Big)+ \\
&&
(z_1/C)^{-m-1}\times\\
&&
\left.
\Big(\Gamma^-_1(z_1)\tau(x+\sqrt\hbar,\qq')\Big)
e^{
  \sum_{\ell>0} \frac{z_1^{(n-2)\ell}}{(n-2)^\ell\ell!}
  (\sqrt\hbar \partial_x +
  h_\ell)(q'_{0,0,\ell}-q''_{0,0,\ell})/\sqrt{\hbar}}
\Big(\Gamma^+_1(z_1)\tau(x-\sqrt\hbar,\qq'')\Big)
\right).
\een
Similarly,
\ben
\widetilde{\Gamma}^{\epsilon^2_1}(\lambda) &= &
\exp\left(-\sum_{\ell=0}^\infty
\frac{\lambda^{\ell+\frac{1}{2}}  }{
\Big(\ell+\frac{1}{2}\Big)\cdots \Big(1+\frac{1}{2}\Big)}\,
(q_{2,1,\ell}+q_{3,1,\ell})/\sqrt\hbar
\right)
e^{\pi\ii\sum_{\ell=0}^\infty
\frac{\lambda^\ell}{\ell!}\,q_{0,0,\ell}/\sqrt{\hbar}
}
\\
&&
\exp\left(\sum_{\ell=0}^\infty
  \Big(0+\frac{1}{2}\Big)\cdots \Big(\ell-1+\frac{1}{2}\Big)
  \lambda^{-\ell-\frac{1}{2}}\,
  \sqrt\hbar (\partial_{2,1,\ell}+\partial_{3,1,\ell})
  \right).
  \een
  Let us make the substitutions: $\lambda=z_2^2/2$ and
  \beq\label{subst-2k}
  t_{2,2\ell+1}=\hbar^{-1/2}\,
  \frac{q_{2,1,\ell}+q_{3,1,\ell}  }{
    \sqrt{2}(2\ell+1)!!}.
  \eeq
  Then, using also \eqref{subst-3k}, we get
  \ben
  \widetilde{\Gamma}^{\pm\epsilon^2_1}(\lambda) =
  e^{\pm\pi\ii\sum_{\ell=0}^\infty
\frac{\lambda^\ell}{\ell!}\,q_{0,0,\ell}/\sqrt{\hbar}
}
\Gamma_2(\mp z_2)
\een
and
\ben
\widetilde{b}_{\epsilon^2_1}(\lambda)
\frac{d\lambda}{\lambda} = -\frac{dz_2}{2z_2}.
\een
We have
\ben
\widetilde{\Gamma}^\#\widetilde{\Gamma}^{\pm\epsilon^2_1}\tau =
\Gamma_2(\mp z_2)\tau(x,\qq)
e^{\pm\sum_{\ell>0}\frac{\lambda^\ell}{\ell!}\,\partial_x \,
  q_{0,0,\ell}}
(-1)^{\pm(q_{0,0,0}+x)/\sqrt\hbar}
\een
and
\ben
\widetilde{\Gamma}\widetilde{\Gamma}^{\pm\epsilon^2_1}\tau =
(-1)^{\pm(q_{0,0,0}+x)/\sqrt\hbar}
e^{\mp\sum_{\ell>0}\frac{\lambda^\ell}{\ell!}\,\partial_x \,
  q_{0,0,\ell}}
\Gamma_2(\mp z_2)\tau(x,\qq).
\een
The residue in the $\lambda$-plane
\ben
\operatorname{Res}_{\lambda=\infty}
\frac{\lambda^r}{r!}\, \widetilde{b}_{\epsilon^2_1}(\lambda)\,
\frac{d\lambda}{\lambda}
\Big(
\widetilde{\Gamma}^\#\otimes \widetilde{\Gamma}
\Big)\Big(
\widetilde{\Gamma}^{\epsilon^2_1}(\lambda)\otimes
\widetilde{\Gamma}^{-\epsilon^2_1}(\lambda) +
\widetilde{\Gamma}^{-\epsilon^2_1}(\lambda)\otimes
\widetilde{\Gamma}^{\epsilon^2_1}(\lambda)
\Big)(\tau\otimes\tau)
\een
turns into the following residue in the $z_2$-plane
\ben
-\operatorname{Res_{z_2=\infty}}
\frac{dz_2}{2z_2}\, \frac{z_2^{2r}}{(2r)!!} \,  (-1)^m
\Big(\Gamma_2(z_2)\tau(x,\qq')\Big)
e^{\sum_{\ell>0}\frac{z_2^{2\ell}}{(2\ell)!!}\,\partial_x \,
  (q'_{0,0,\ell}-q''_{0,0,\ell})}
\Big(\Gamma_2(-z_2)\tau(x,\qq'')\Big).
\een
Finally
\ben
\widetilde{\Gamma}^{\epsilon^3_1}(\lambda) &= &
\exp\left(-\sum_{\ell=0}^\infty
\frac{\lambda^{\ell+\frac{1}{2}}  }{
\Big(\ell+\frac{1}{2}\Big)\cdots \Big(1+\frac{1}{2}\Big)}\,
(q_{2,1,\ell}-q_{3,1,\ell})/\sqrt\hbar
\right)
\\
&&
\exp\left(\sum_{\ell=0}^\infty
  \Big(0+\frac{1}{2}\Big)\cdots \Big(\ell-1+\frac{1}{2}\Big)
  \lambda^{-\ell-\frac{1}{2}}\,
  \sqrt\hbar (\partial_{2,1,\ell}-\partial_{3,1,\ell})
  \right).
  \een
  Let us make the substitutions: $\lambda=z_3^2/2$ and
  \beq\label{subst-3k}
  t_{3,2\ell+1}=\hbar^{-1/2}\,
  \frac{q_{2,1,\ell}-q_{3,1,\ell}  }{
    \sqrt{2}(2\ell+1)!!}.
  \eeq
  Then, using also \eqref{subst-2k}, we get
$
\widetilde{\Gamma}^{\pm\epsilon^3_1}(\lambda) =
\Gamma_3(\mp z_3)
$
and
\ben
\widetilde{b}_{\epsilon^3_1}(\lambda)
\frac{d\lambda}{\lambda} = \frac{dz_2}{2z_2}.
\een
We have
\ben
\widetilde{\Gamma}^\#\widetilde{\Gamma}^{\pm\epsilon^3_1}\tau =
\Gamma_3(\mp z_3)\tau(x,\qq)
e^{\pm\sum_{\ell>0}\frac{\lambda^\ell}{\ell!}\,\partial_x \,
  q_{0,0,\ell}}
\een
and
\ben
\widetilde{\Gamma}\widetilde{\Gamma}^{\pm\epsilon^3_1}\tau =
e^{\mp\sum_{\ell>0}\frac{\lambda^\ell}{\ell!}\,\partial_x \,
  q_{0,0,\ell}}
\Gamma_3(\mp z_3)\tau(x,\qq).
\een
The residue in the $\lambda$-plane
\ben
\operatorname{Res}_{\lambda=\infty}
\frac{\lambda^r}{r!}\, \widetilde{b}_{\epsilon^3_1}(\lambda)\, 
\frac{d\lambda}{\lambda}
\Big(
\widetilde{\Gamma}^\#\otimes \widetilde{\Gamma}
\Big)\Big(
\widetilde{\Gamma}^{\epsilon^3_1}(\lambda)\otimes
\widetilde{\Gamma}^{-\epsilon^3_1}(\lambda) +
\widetilde{\Gamma}^{-\epsilon^3_1}(\lambda)\otimes
\widetilde{\Gamma}^{\epsilon^3_1}(\lambda)
\Big)(\tau\otimes\tau)
\een
turns into the following residue in the $z_3$-plane
\ben
\operatorname{Res_{z_3=\infty}}
\frac{dz_3}{2z_3}\, \frac{z_3^{2r}}{(2r)!!} \,
\Big(\Gamma_3(z_3)\tau(x,\qq')\Big)
e^{\sum_{\ell>0}\frac{z_3^{2\ell}}{(2\ell)!!}\,\partial_x \,
  (q'_{0,0,\ell}-q''_{0,0,\ell})}
\Big(\Gamma_3(-z_3)\tau(x,\qq'')\Big).
\een
\subsection{Wave functions}
Let $\O_\epsilon(\CC):=\O(\CC)[\![\epsilon]\!]$ denote the ring of formal power series in
$\epsilon$ whose coefficients are holomorphic functions on $\CC$. Let
us denote by $x$ the standard coordinate function on $\CC$. We will be
interested in formal functions of the form 
$\tau(x,\tt)=e^{\F (\epsilon, x, \tt)\epsilon^{-1}}$, where
$\F(\epsilon, x, \tt)$  is a formal power series in $\tt$ with
coefficients in $\O_\epsilon(\CC)$ satisfying the following condition:
the coefficient in front of $\epsilon^0$ in $\F(\epsilon, x, \tt)$ is
at most linear in $\tt$. For example, if $\D(\hbar,\qq)$ is the total
descendent potential of $\PP^1_{n-2,2,2}$, then the substitution
$(\hbar,\qq)\mapsto (\epsilon, x,\tt)$ defined by 
\ben
\sqrt{\hbar}:=\epsilon,\quad
q_{0,0,0}:=x,\quad
q_{0,0,k}:= t_{0,k} \epsilon -\delta_{k,1}\ (k>0)
\een
and the linear change of variables \eqref{subst-1k1}--\eqref{subst-1k2},
\eqref{subst-2k}, and \eqref{subst-3k}, identify the total descendent
potential with a formal function of the above type. In order to prove
this, we need only to recall the dimension formula for the virtual
fundamental cycle, that is, if the Gromov--Witten invariant
\ben
\langle \phi_{i_1,p_1} \psi^{k_1},\dots,\phi_{i_r,p_r}\psi^{k_r}\rangle_{g,r,d}
\een
is not 0, then
\ben
\sum_{s=1}^r \left( k_s + \tfrac{p_s}{a_{i_s}}\right) = 2g-2+r +\frac{d}{n-2}.
\een
Let us point out that, using the above formula, we can prove also that the
coefficient in front of each monomial in $\tt$ and $\epsilon$ in $\log
\D(\hbar,\qq)$ is a polynomial in $x$ and $Q$.

Let us introduce the following functions:
\ben
\Psi^+_1(x,\tt,z) :=  \psi^+_1(x,\tt,z)e^{\xi_1(\tt,z)}
(-z/C)^{\tfrac{x}{\epsilon}-\tfrac{1}{2}},\quad
\Psi^-_1(x,\tt,z) :=  (-z/C)^{-\tfrac{x}{\epsilon} -\frac{1}{2}}
e^{-\xi_1(\tt,z)}\psi^-_1(x,\tt,z),
\een
where
\ben
\xi_1(\tt,z) & :=& \sum_{k=1}^\infty \Big(
t_{1,k} z^k + t_{0,k}\,\frac{z^{(n-2)k}  }{(n-2)^kk!}
(\epsilon \partial_x-h_k)\Big),
\een
\ben
\psi_1^\pm(x,\tt,z) :=
\frac{e^{\mp\sum_{k=1}^\infty \frac{z^{-k}}{k} \,\partial_{t_{1,k}}}
  \tau(x\mp\epsilon,\tt)}{\tau(x,\tt)} =:\sum_{k=0}^\infty
\psi^{\pm}_{1,k}(x,\tt) z^{-k},
\een
and
\ben
\Psi^+_a(x,\tt,z)  :=  \psi^+_a(x,\tt,z) e^{\xi_a(\tt,z)},\quad
\Psi^-_a(x,\tt,z) :=  e^{-\xi_a(\tt,z)}\psi^-_a(x,\tt,z) \quad
(a=2,3), 
\een
where
\ben
\xi_a(\tt,z) & :=& \sum_{k=1}^\infty \Big(
t_{a,2k-1} z^{2k-1} + t_{0,k}\, \frac{ z^{2k} }{2^k k!}\,\epsilon\partial_x \Big),
\een
\ben
\psi^\pm_a(x,\tt,z):=\frac{
  e^{\mp 2 \, \sum_{k=0}^\infty
    \frac{z^{-2k-1}}{2k+1}\,\partial_{t_{a,2k+1}}}
  \tau(x,\tt)}{
  \tau(x,\tt)} =:
\sum_{k=0}^\infty \psi^{\pm}_{a,k}(x,\tt) z^{-k}.
\een
The functions $\Psi^\pm_i(x,\tt,z)$ ($1\leq i\leq 3$) will be called {\em wave
  functions} of the Extended D-Toda Hierarchy if they satisfy the
following bilinear equations: 
\ben
\operatorname{Res}_{z=0} \frac{z^{(n-2)r}}{(n-2)^rr!}
\frac{dz}{z}\Big(
\Psi^+_1(x,\tt',z)\Psi^-_1(x+m\epsilon,\tt'',z) +
(\Psi^+_1(x+m\epsilon ,\tt'',z)\Psi^-_1(x,\tt',z))^\#\Big)=\\
\operatorname{Res}_{z=0} \frac{z^{2r}}{2^r r!}
\frac{dz}{2z}\Big(
\Psi^+_2(x,\tt',z)\Psi^-_2(x+m \epsilon,\tt'',z) -(-1)^{m}
\Psi^+_3(x,\tt',z)\Psi^-_3(x+m \epsilon,\tt'',z)
\Big),
\een
where $r\geq 0$ and $m$ are arbitrary integers.

Suppose that $\D(\hbar,\qq)\in
\CC_\hbar[\![q_0,q_1+\mathbf{1},q_2,\dots ]\!]$ is a formal power
series, such that, under the substitution $(\hbar,\qq)\mapsto
(\epsilon,x,\tt)$ from 
above, $\D(\hbar,\qq)$ becomes a formal function in $\epsilon,x,$ and
$\tt$ for which the definition of $\Psi^\pm_i(x,\tt,z)$ ($1\leq i\leq
3$) makes sense. Then one can check that
the corresponding functions $\Psi^\pm_i(x,\tt,z)$ ($1\leq i\leq 3$)
are wave functions of the Extended D-Toda Hierarchy if and only if 
$\D(\hbar,\tt)$ satisfies the HQEs
\eqref{hqe}.

\section{Phase factors}

The proof of Theorem \ref{t1} follows the method developed by Givental
in \cite{G3}. It relies on Givental's higher genus
reconstruction formula (see \cite{G1}) proved by Teleman \cite{Te}. 
There is a part of Givental's argument in \cite{G3} which was hard to
generalize. The difficulty however was offset in \cite{Mi2} in the
settings of singularity theory and after a small modification in the
current settings as well (see also \cite{MST}). The goal of this section is
to prove the results that we need in order to make
Givental's argument from \cite{G3} work.

\subsection{Proof of Lemma \ref{le:ip}}\label{sec:ip}
The calibration $S(t,Q,z)$ satisfies the so-called {\em divisor
  equation}
\ben
\partial_{t_{0,1}} S(t,Q,z) = Q\partial_Q S(t,Q,z) + S(t,Q,z) \,
z^{-1} P\cup,
\een
where $P\cup$ is the linear operator defined by topological cup
product by the cohomology class $\phi_{0,1}=P$. In order to avoid
cumbersome notation, let us agree to denote by $P$ the linear operator
$P\cup$ and by $\phi_{0,1}$ the cohomology class $P$. Using the
commutation relation $\theta P=P (\theta-1)$ and \eqref{transl-invar}, we get
\ben
\widetilde{I}^{(m)}(\lambda) P= P
\partial_\lambda \widetilde{I}^{(m)}(\lambda). 
\een
The above formula implies that
\ben
I^{(m)}(t,\lambda) e^{-P\log Q} a =
S(t,Q,-\partial_\lambda^{-1}) e^{-P\log Q \partial_\lambda}
\widetilde{I}^{(m)}(\lambda) a.
\een
On the other hand, according to the divisor equation, we have
$Q\partial_Q\left(S(t,Q,z) e^{P\log Q/z} \right) = z^{-1}P\bullet
S(t,Q,z) e^{P\log Q/z}$. Therefore,
\ben
Q\partial_Q\ \left(
I^{(m)}(t,\lambda) e^{-P\log Q} a\right) =
-P\bullet I^{(m+1)}(t,\lambda) e^{-P\log Q} a.
\een
The above differential equation and
$(\lambda-E\bullet)I^{(m+1)} =
(\theta-m-\tfrac{1}{2}) I^{(m)}$ imply that the pairing
\ben
( e^{-P\log Q} a | e^{-P\log Q} b ) = \left(
  I^{(0)}(t,\lambda) e^{-P\log Q} a ,
  (\lambda-E\bullet)
  I^{(0)}(t,\lambda) e^{-P\log Q} b
  \right)
  \een
  is independent of $Q$. On the other, hand the operators $S_k(t,Q)$
  vanish at $t=Q=0$. Therefore, the only contribution to the
  intersection pairing is given by
  \beq\label{ip-regul}
  ( e^{-P\log Q} a | e^{-P\log Q} b ) =
  \left(
  \widetilde{I}^{(0)}(\lambda) a ,
  (\lambda-\rho)
  \widetilde{I}^{(0)}(\lambda) b
  \right).
  \eeq
  On the other hand, since $\rho$ is a multiplication by a top degree
  class and $\theta(P) = -\tfrac{1}{2}\, P$, we have
  $\left(\theta +\tfrac{1}{2}\right) \rho = 0$. This relation allows us
  to get the following simple formula for the calibrated periods (for
  $m<0$) 
  \beq\label{calibr-period}
  \widetilde{I}^{(m)} (\lambda) =
  \frac{
    \lambda^{\theta-m-\tfrac{1}{2}}  } {
    \Gamma(\theta-m+\tfrac{1}{2} ) } +
  \frac{\lambda^{-m-1}} {\Gamma(-m)}\, \left(
    \log\lambda -\psi(-m)\right) \, \rho,
  \eeq
  where $\psi(x):=\tfrac{\Gamma'(x)}{\Gamma(x)}$ is the di-gamma
  function. From here, since
  $\widetilde{I}^{(0)}=\partial_\lambda \widetilde{I}^{(-1)}$, we get
  \ben
  \widetilde{I}^{(0)} (\lambda) =
  \frac{
    \lambda^{\theta-\tfrac{1}{2}}  } {
    \Gamma(\theta+\tfrac{1}{2} ) } +\frac{1}{\lambda} \, \rho.
  \een
  Substituting the above formula in \eqref{ip-regul}, after some
  straightforward computation, we get
  \ben
  ( e^{-P\log Q} a | e^{-P\log Q} b ) =\frac{1}{2\pi}
  \left( (e^{\pi\sqrt{-1} \theta } +
     e^{-\pi\sqrt{-1} \theta } + 2\pi \rho) a , b\right).
    \een
   In order to complete the proof of Lemma \ref{le:ip}, it remains
   only to notice the 
   following two facts: the RHS of the above formula is invariant
   under the action $(a,b)\mapsto (e^{P\log Q}a, e^{P\log Q}b)$ and 
   \ben
   e^{\pi\sqrt{-1} \theta } +
   e^{-\pi\sqrt{-1} \theta } + 2\pi \rho =
   e^{\pi\sqrt{-1} \theta } e^{\pi\sqrt{-1} \rho } +
   e^{\pi\sqrt{-1} \rho } e^{-\pi\sqrt{-1} \theta } ,
   \een
   where we used the relations $\rho^2=0$ and $\theta \rho
   =-\tfrac{1}{2} \rho$.

   \subsection{Euler pairing}\label{sec:eup}
   
Let us denote the imaginary unit by $\ii:=\sqrt{-1}$. 
In what follows, we will have to use quite frequently the following
formulas for the Euler pairing:
\ben
\langle a,b\rangle:= \chi(a^\vee\otimes b),\quad a,b\in K^0(\PP^1_{n-2,2,2}).
\een
If $1\leq i,j \leq n-3$, then
\ben
\langle \epsilon^1_i,\epsilon^1_j\rangle = 
\begin{cases}
\phantom{-}\frac{1}{2} & \mbox{if}\quad i\leq j, \\
-\frac{1}{2} & \mbox{if} \quad i>j.
\end{cases}
\een
If $1\leq i\leq n-3$, then
\ben
\langle \epsilon^1_i,\epsilon^2_1\rangle = \frac{1}{2},\quad
\langle \epsilon^2_1,\epsilon^1_i\rangle = -\frac{1}{2},
\quad
\langle \epsilon^1_i,\epsilon^3_1\rangle = 
\langle \epsilon^3_1,\epsilon^1_i\rangle =0.
\een
Finally
\ben
\langle \epsilon^2_1,\epsilon^2_1\rangle=
\langle \epsilon^3_1,\epsilon^3_1\rangle=
\frac{1}{2},\quad
\langle \epsilon^2_1,\epsilon^3_1\rangle=
\langle \epsilon^3_1,\epsilon^2_1\rangle=0.
\een
Let us sketch the proof of the first formula. The computations in the
remaining ones are similar and much shorter. Note that
\ben
\Psi(\epsilon^1_i) = 1+
\frac{\pi\ii(2i+n-2)-\gamma}{n-2}\, P +
\sum_{p=1}^{n-3} \Gamma\Big(1-\frac{p}{n-2}\Big)\eta_1^{-p i}
\phi_{1,p},
\een
and
\ben
e^{\pi\ii\theta}e^{\pi\ii\rho}\Psi(\epsilon^1_j) = 
\ii + 
\frac{\pi (2j+n-1) +\ii\gamma }{n-2}\,  P + 
\sum_{p=1}^{n-3} \Gamma\Big(1-\frac{p}{n-2}\Big)\eta_1^{-p(j+1/2)}
\ii\phi_{1,p}
\een
where $\eta_1=e^{2\pi\ii/(n-2)}$. Recalling formula
\eqref{euler-pairing}, we get 
\ben
\langle \epsilon^1_i,\epsilon^1_j\rangle = 
\frac{1}{n-2}\Big(
j-i+\frac{1}{2} - \sum_{p=1}^{n-3} \frac{\eta_1^{(i-j)p}}{\eta_1^p-1}
\Big).
\een
Suppose that $i-j>0$ then the sum 
\ben
\sum_{p=1}^{n-3} \frac{\eta_1^{(i-j)p}}{\eta_1^p-1} = 
\sum_{p=1}^{n-3} \Big(\frac{1}{\eta_1^p-1} +
1+\eta_1^p+\cdots+\eta_1^{p(i-j-1)}\Big) = j-i+\frac{1}{2} +\frac{n-2}{2}.
\een
If $j\geq i$ then since the above sum is invariant under $p\mapsto -p$ we get
\ben
\sum_{p=1}^{n-3} \frac{\eta_1^{(i-j)p}}{\eta_1^p-1}=
\sum_{p=1}^{n-3} \frac{\eta_1^{-(i-j)p}}{\eta_1^{-p}-1} =
-\sum_{p=1}^{n-3} \frac{\eta_1^{(1-i+j)p}}{\eta_1^{p}-1}=j-i+\frac{1}{2}-\frac{n-2}{2}.
\een
The formulas for $\langle\epsilon^1_i,\epsilon^1_j\rangle$ follow. 

Let us fix an integer $\kappa$ divisible by the order of the semi-simple part of
$\sigma$. Using Proposition \ref{pr:cal-periods} we get that 
$\kappa=2(n-2)$ is such a number. If $\beta\in H$, then let us define
\ben
\beta_0  := 
\frac{1}{\kappa}(\beta+\sigma(\beta)+\cdots+\sigma^{\kappa-1}(\beta)),\quad
\beta_{\rm tw} :=  \beta-\beta_0.
\een
\begin{lemma}\label{le:a-euler-b0}
The following formulas hold
\ben
\langle \alpha,\beta_0\rangle=
\operatorname{rk}(\alpha)\operatorname{deg}(\beta)-
\operatorname{rk}(\beta)\operatorname{deg}(\alpha) + 
\operatorname{rk}(\alpha)\operatorname{rk}(\beta)
\een
and
\ben
\langle \beta_0,\alpha\rangle=
\operatorname{rk}(\beta)\operatorname{deg}(\alpha)-
\operatorname{rk}(\alpha)\operatorname{deg}(\beta) + 
\operatorname{rk}(\alpha)\operatorname{rk}(\beta)\Big(\frac{1}{n-2}-1\Big).
\een
\end{lemma}
\proof
Using that $\chi_{j,p}(\sigma(E)) = \eta_j^{-p}\chi_{j,p}(E)$, where
$\eta_j=e^{2\pi\ii/a_j}$ we get that
\ben
e^{\pi\ii\theta}e^{\pi\ii\rho}\Psi(\beta_0)=\ii\operatorname{rk}(\beta)+\Big(
\operatorname{rk}(\beta)
\frac{\pi+\ii\gamma}{n-2} +
2\pi\operatorname{deg}(\beta_0)\Big) P.
\een
Since
\ben
\Psi(\alpha)=
\operatorname{rk}(\alpha)+\Big(
-\operatorname{rk}(\alpha)
\frac{\gamma}{n-2} +
2\pi\ii\operatorname{deg}(\alpha)\Big) P
+\cdots,
\een
where the dots involve cohomology classes supported on the twisted
sectors. Recalling again formula \eqref{euler-pairing} we get
\ben
\langle\alpha,\beta_0\rangle =
\operatorname{rk}(\alpha)\operatorname{deg}(\beta_0)-
\operatorname{rk}(\beta)\operatorname{deg}(\alpha) + 
\frac{1}{\kappa}\, \operatorname{rk}(\alpha)\operatorname{rk}(\beta)
\een
Since $\sigma$ is multiplication by the tangent bundle
$T\PP^1_{n-2,2,2}$ we get that
\ben
\operatorname{deg}(\beta_0) =
\operatorname{deg}(\beta)+\operatorname{rk}(\beta)\frac{\kappa-1}{\kappa}.
\een
The formula for $\langle\alpha,\beta_0\rangle$ follows. The second
formula is proved in a similar way, so we omit the argument.
\qed
\begin{lemma}\label{le:ab_tw}
  The following formula holds
  \ben
  (\alpha|\sigma^s(\beta_{\rm tw})) = \sum_{j=1}^3\sum_{p=1}^{a_j-1}
  \frac{1}{a_j}\eta_j^{ps}\chi_{j,p}(\alpha)\chi_{j,a_j-p}(\beta).
  \een
\end{lemma}
\proof
Note that
\ben
\operatorname{rk}(\beta_{\rm tw}) = 0,\quad
\operatorname{deg}(\beta_{\rm tw}) = -
\operatorname{rk}(\beta)(1-1/\kappa),
\een
and
\ben
\Psi(\sigma^s\beta_{\rm tw}) = -2\pi\ii
\operatorname{rk}(\beta)(1-1/\kappa)P+
\sum_{j,p} \Gamma\Big(1-\frac{p}{a_j}\Big)\chi_{j,p}(\beta)\eta_j^{-ps}\phi_{j,p}.
\een
Recalling again formula \eqref{euler-pairing} we get that the pairing
$(\alpha|\sigma^s(\beta_{\rm tw}))$ is given by
\ben
\frac{1}{2\pi}(\Psi(\alpha),(e^{-\pi\ii\theta}+
e^{\pi\ii\theta})\Psi(\sigma^s\beta_{\rm tw}))=
\sum_{j,p}(\Psi(\alpha),\phi_{j,p})
\Gamma(p/a_j)^{-1}\chi_{j,p}(\beta)\eta_j^{-ps}.
\een
It remains only to use that
\ben
(\Psi(\alpha),\phi_{j,p}) = \frac{1}{a_j} \Gamma(p/a_j)
\chi_{j,a_j-p}(\alpha).
\qed
\een
\subsection{Phase factors}\label{sec:pf}
By definition phase factors are scalar functions that arise when we
compose two vertex operators. In our case, the phase factors can be
expressed in terms of Givental's symplectic form as follows: 
\ben
\widetilde{\Gamma}^\alpha(\lambda_1)
\widetilde{\Gamma}^\beta (\lambda_2)  =
e^{\Omega(\widetilde{\ff}_\alpha^+(\lambda_1,z),
  \widetilde{\ff}_\beta^-(\lambda_2,z) )}\ 
\normord
\widetilde{\Gamma}^\alpha(\lambda_1)
\widetilde{\Gamma}^\beta (\lambda_2)
\normord,
\een
and
\ben
{\Gamma}^\alpha(t,\lambda_1)
{\Gamma}^\beta (t,\lambda_2)  =
e^{\Omega({\ff}_\alpha^+(t,\lambda_1,z),
  {\ff}_\beta^-(t,\lambda_2,z) )} \   
\normord
\Gamma^\alpha(t,\lambda_1)
\Gamma^\beta (t,\lambda_2)
\normord,
\een
where the normal ordering $\normord\ \normord$ means that all
differentiations should be moved to the right.
\begin{proposition}\label{ope-cal-periods}
If $\alpha,\beta\in K^0(\PP^1_{n-2,2,2})$ then the symplectic pairing
\ben
&&
\Omega(\widetilde{\ff}_\alpha^+(\lambda_1,z),
\widetilde{\ff}_\beta^-(\lambda_2,z) )= -2\pi\ii\,
\operatorname{rk}(\alpha)\operatorname{deg}(\beta) +\\
&&
+\log \Big(
(Q\lambda_2^{-1/(n-2)})^{\operatorname{rk}(\alpha)
  \operatorname{rk}(\beta)}
\prod_{s=1}^\kappa
\Big(1-\eta^{-s}(\lambda_2/\lambda_1)^{1/\kappa}\Big)^
{(\alpha|\sigma^s\beta)}\Big),
\een
where $\eta=e^{2\pi\ii/\kappa}$ and the RHS should be expanded into a
Laurent series in $\lambda_1^{-1}$ in the region
$|\lambda_1|>|\lambda_2|$. 
\end{proposition}
\proof
By definition the symplectic pairing is
\ben
\sum_{m=0}^\infty
(-1)^{m+1}(
\widetilde{I}^{(m)}_\alpha(\lambda_1), 
\widetilde{I}^{(-m-1)}_\beta(\lambda_2)).
\een
Substituting the explicit formulas for the calibrated periods from
Proposition \ref{pr:cal-periods}, part a) we get
\beqa\nonumber
-2\pi\ii \operatorname{rk}(\alpha)\operatorname{deg}(\beta) -
\frac{1}{n-2}\operatorname{rk}(\alpha)\operatorname{rk}(\beta)
\Big(\log \lambda_2-C_0+
\sum_{\ell=1}^\infty \frac{(\lambda_2/\lambda_1)^\ell}{\ell}
\Big)\\
\label{2ndline}
-\sum_{j=1}^3\sum_{p=1}^{a_j-1}\sum_{m=0}^\infty
\frac{(\lambda_2/\lambda_1)^{m+p/a_j}}
{p+m a_j}\, \chi_{j,p}(\alpha)\chi_{j,a_j-p}(\beta).
\eeqa
Let us compute the sum over $m$ in \eqref{2ndline}. We will use the
following simple fact: if $r$ is an arbitrary integer, then 
\ben
\frac{1}{\kappa} \sum_{s=1}^\kappa \eta^{(r-l)s} =
\begin{cases}
1 & \mbox{ if  } l=r+m\kappa \mbox{ for some integer } m,\\
0 & \mbox{otherwise}.
\end{cases}
\een
Using the above formula we get the following identity:
\ben
\frac{1}{\kappa} \sum_{l=1}^\infty \sum_{s=1}^\kappa \eta^{rs} 
\frac{(\eta^{-s}(\lambda_2/\lambda_1)^{1/\kappa})^l}{l} = 
\sum_{m=0}^\infty
\frac{(\lambda_2/\lambda_1)^{m+\tfrac{r}{\kappa}}}{r+m\kappa} .
\een
Let us specialize the above identity to $r=p\kappa/a_j$. Since
$\eta^{rs} = e^{2\pi\ii rs/\kappa} = \eta_j^{ps}$, we get 
\ben
\frac{1}{a_j} \sum_{l=1}^\infty \sum_{s=1}^\kappa \eta_j^{ps} 
\frac{(\eta^{-s}(\lambda_2/\lambda_1)^{1/\kappa})^l}{l} = 
\sum_{m=0}^\infty
\frac{(\lambda_2/\lambda_1)^{m+\tfrac{p}{a_j} } }{p+m a_j} .
\een
Using the above formula and Lemma \ref{le:ab_tw}, we get that the
expression in \eqref{2ndline} is equal to
\ben
-\sum_{\ell=1}^\infty\sum_{s=1}^\kappa
\frac{1}{\ell} (\eta^{-s} (\lambda_2/\lambda_1)^{1/\kappa})^\ell
(\alpha|\sigma^s\beta_{\rm tw}) =
\log \prod_{s=1}^\kappa
\Big(1-\eta^{-s} (\lambda_2/\lambda_1)^{1/\kappa}\Big)^{
(\alpha|\sigma^s\beta_{\rm tw})}.
\een
Recalling Lemma \ref{le:a-euler-b0}, we get
\ben
(\alpha|\beta_0) =
\frac{1}{n-2}\operatorname{rk}(\alpha) \operatorname{rk}(\beta).
\een
Note that this formula implies also that
$(\alpha|\sigma^s\beta_0)=(\alpha|\beta_0)$. 
Finally, since $C_0=(n-2)\log Q$ and
\ben
-\sum_{\ell=1}^\infty \frac{(\lambda_2/\lambda_1)^\ell}{\ell} =
\log (1-\lambda_2/\lambda_1) = \log
\prod_{s=1}^\kappa
\Big(1-\eta^{-s} (\lambda_2/\lambda_1)^{1/\kappa}\Big)
\een
the formula that we want to prove follows.
\qed

\begin{lemma}\label{le:aux-id}
The following formula holds:
\ben
\sum_{s=1}^\kappa \Big(\frac{1}{2}-\frac{s}{\kappa}\Big)
(\alpha|\sigma^s\beta)=
-\langle \alpha,\beta\rangle + 
\operatorname{rk}(\alpha) \operatorname{deg}(\beta)-
\operatorname{rk}(\beta) \operatorname{deg}(\alpha).
\een
\end{lemma}
\proof
Since $\sigma^{-1}$ is multiplication by the dualizing sheaf, we can
write Serre's duality as 
\ben
\langle b,a\rangle =-\langle a,\sigma^{-1}b\rangle.
\een
Therefore, the intersection pairing 
\ben
(a|b) = \langle a, (1-\sigma^{-1})b\rangle,\quad a,b\in K^0(\PP^1_{n-2,2,2}).
\een
Using the above formula let us write the LHS of the identity that we
want to prove as 
\beq\label{aux-id-1}
\sum_{s=1}^\kappa \Big(\frac{1}{2}-\frac{s}{\kappa}\Big)
\langle\alpha, (1-\sigma^{-1})\sigma^s\beta\rangle =
\langle\alpha,-\beta -\frac{1}{2}(\sigma^k-1)\beta +\beta_0\rangle.
\eeq
Here we wrote $(1-\sigma^{-1}) \sigma^s=\sigma^s-\sigma^{s-1}$,
split the sum over $s$ into two sum. After shifting  the summation
index in the second sum $s\mapsto s+1$, most of the terms cancel out
and what is left is the RHS of \eqref{aux-id-1}. Note that 
\ben
\langle \alpha, (\sigma^\kappa-1)\beta\rangle = 
\kappa \langle \alpha, (\sigma-1)\beta_0\rangle =
\kappa (\sigma^{-1}\alpha|\beta_0) = 
2\operatorname{rk}(\alpha)\operatorname{rk}(\beta),
\een
where we used Lemma \ref{le:a-euler-b0} for the last
equality. Therefore \eqref{aux-id-1} becomes 
\ben
-\langle\alpha,\beta\rangle -
\operatorname{rk}(\alpha)\operatorname{rk}(\beta) +
\langle\alpha, \beta_0\rangle.
\een
Recalling Lemma \ref{le:a-euler-b0} one
more time completes the proof. 
\qed

\medskip
Let us discuss the analytic properties of the phase factors. The
symplectic pairing 
$\Omega(\widetilde{\ff}^+_\alpha(\lambda_1,z),
\widetilde{\ff}^-_\beta(\lambda_2,z) )$ is a convergent Laurent
series that extends analytically to a multivalued analytic function
$\widetilde{\Omega}_{\alpha,\beta}(\lambda_1,\lambda_2)$ in the domain
\ben
0<|\lambda_1-\lambda_2|<\operatorname{min}(|\lambda_1|,|\lambda_2|).
\een
Using Proposition \ref{ope-cal-periods} and Lemma \ref{le:aux-id} we
get that  
\beq\label{calibr-locality}
\widetilde{\Omega}_{\alpha,\beta}(\lambda_1,\lambda_2)-
\widetilde{\Omega}_{\beta,\alpha}(\lambda_2,\lambda_1) =
-2\pi\ii \Big(
\langle \alpha,\beta\rangle + \ell (\alpha|\beta)\Big)
\eeq
for some integer $\ell\in\ZZ$ that depends on how we fix the branches
of the functions on the LHS.

The symplectic pairing
\beq\label{phase-factor}
\Omega(\ff_\alpha^+(t,\lambda_1,z), \ff_\beta^-(t,\lambda_2,z)) =
\sum_{m=0}^\infty (-1)^{m+1}
(I^{(m)}_\alpha(t,\lambda_1), I^{(-m-1)}_\beta(t,\lambda_2))
\eeq
will be interpreted via its Laurent series expansion in
$\lambda_1^{-1}.$ The radius of convergence will be determined
below. Let us express \eqref{phase-factor} in terms of
$\Omega(\widetilde{\ff}_\alpha^+(\lambda_1,z),
\widetilde{\ff}_\beta^-(\lambda_2,z)) $ and the calibration operator
$S(t,Q,z)$. We have the following conjugation formula
(see \cite{G3}, formula (17)):
\beq\label{S-conj}
\widehat{S} e^{\widehat{\ff}} \widehat{S}^{-1} =
e^{W(\ff^+,\ff^+)/2} e^{(S\ff)\sphat},
\eeq
where $W(t,\qq,\qq)=\sum_{k,\ell=0}^\infty (W_{k,\ell}(t)q_\ell,q_k)$ is the quadratic
form defined by
\ben
\sum_{k,\ell=0}^\infty W_{k,\ell}(t) w^{-k} z^{-\ell} =
\frac{S^*(t,Q,w)S(t,Q,z)-1}{z^{-1}+w^{-1}},
\een
where ${}^*$ denotes transposition with respect to the Poincare
pairing, and if $\ff=\sum_{k=0}^\infty f_k z^k$, then
$W(t,\ff,\ff):=\sum_{k,l= 0}^\infty (W_{k,l}(t) f_l,f_k). $ Put 
\beq\label{Wab}
W_{\alpha,\beta}(t,\lambda_1,\lambda_2):=
W(t,
\widetilde{\ff}^+_\alpha(\lambda_1,z),
\widetilde{\ff}^+_\beta(\lambda_2,z)).
\eeq
Since
$\ff_\alpha(t,\lambda,z)=S(t,Q,z)\widetilde{\ff}_\alpha(\lambda,z)$,
the conjugation formula \eqref{S-conj} and the the definition of the
phase factors, i.e., the fact that the symplectic pairings arise as a
composition of vertex operators yield the following relation:
\beq\label{phase-Wab}
\Omega(\ff_\alpha^+(t,\lambda_1,z), \ff_\beta^-(t,\lambda_2,z)) =
\Omega(
\widetilde{\ff}_\alpha^+(\lambda_1,z),
\widetilde{\ff}_\beta^-(\lambda_2,z) ) +
W_{\alpha,\beta}(t,\lambda_1,\lambda_2).
\eeq
On the other hand, since $S(t,Q,z)$ is a solution to the Dubrovin's
connection with respect to $t$, we get that the de Rham differential
(with respect to $t\in H$)
\ben
d W_{\alpha,\beta}(t,\lambda_1,\lambda_2) =
I^{(0)}_\alpha(t,\lambda_1)\bullet
I^{(0)}_\beta (t,\lambda_2),
\een
where the RHS is identified with an element in $T_t^*H=H^*$ via the
Poincare pairing. Using the divisor equation, we get that  $\lim_{t\to -\infty}
S(t,Q,z)e^{-t_{0,1} P/z}= 1$ where we define $t\to -\infty$ to be the
limit $t_{i,p}\to 0$ ($(i,p)\neq (0,1)$) and
$\operatorname{Re}(t_{0,1})\to -\infty$. Although the limit of
$W_{\alpha,\beta}(t,\lambda_1,\lambda_2)$ does 
not exists, it is not hard to analyse the singular term, i.e., if we
subtract $\tfrac{1}{n-2} t_{0,1}
\operatorname{rk}(\alpha)\operatorname{rk}(\beta)$, then the limit
will exists and it is 0. Therefore, we have the following integral
formula:
\beq\label{Wab-integral}
W_{\alpha,\beta}(t,\lambda_1,\lambda_2) =
\tfrac{1}{n-2} t_{0,1}
\operatorname{rk}(\alpha)\operatorname{rk}(\beta) +
\int_{-\infty}^t \Big(
I^{(0)}_\alpha(t',\lambda_1)\bullet
I^{(0)}_\beta (t',\lambda_2) - \tfrac{1}{n-2} dt_{0,1}'
\operatorname{rk}(\alpha)\operatorname{rk}(\beta)
\Big).
\eeq
Put $r(t):=\operatorname{max}_i |u_i(t)|$, where $\{u_i(t)\}$ is the
set of eigenvalues of $E\bullet_t$. We claim that the infinite series
\eqref{phase-factor} is
convergent for all $(t,\lambda_1,\lambda_2)$ in the domain 
\ben
D^+_\infty:=\{
(t,\lambda_1,\lambda_2)\in H\times \CC^2\quad |\quad 
|\lambda_1-\lambda_2|<|\lambda_1|-r(t),\quad |\lambda_2|<|\lambda_1|\}.
\een
Indeed, the period vectors $I^{(m)}(t,\lambda)$ are solutions to a
Fuchsian connection that has poles at $\lambda=u_i(t)$. Therefore,
their Laurent series expansions at $\lambda=\infty$ are convergent for
$|\lambda|>r(t)$. The symplectic pairing
$\Omega(\widetilde{\ff}_\alpha(\lambda_1,z),
\widetilde{\ff}_\beta(\lambda_2,z))$ is convergent in the domain
$|\lambda_1|>|\lambda_2|$ (see Proposition
\ref{ope-cal-periods}). Combining these observations with formulas
\eqref{phase-Wab} and \eqref{Wab-integral}, we get that the Laurent
series expansion at $\lambda_1=\infty$ of \eqref{phase-factor} is convergent for all
$(t,\lambda_1,\lambda_2)$ satisfying
$|\lambda_1|>|\lambda_2|$ and $|\lambda_1|>r(t) $. Clearly, the domain
$D^+_\infty$ satisfies these 
inequalities. The reason why we imposed also the inequality 
$|\lambda_1-\lambda_2|<|\lambda_1|-r(t)$ is to guarantee that the
straight segment from $\lambda_1$ to $\lambda_2$ is outside the disk
with center $0$ and radius $r(t)$, which allows us to make the 
value of the symplectic pairing \eqref{phase-factor} depends only on the choice of a
reference path from $(t^\circ,\lambda^\circ)$ to $(t,\lambda_1)$. Such
a path determines the value of $I^{(m)}_\alpha(t,\lambda_1)$, while
the composition of the reference path and the straight segment from
$(t,\lambda_1)$ to $(t,\lambda_2)$ provides a reference path that
fixes the value of the periods $I^{(-m-1)}_\beta(t,\lambda_2)$. Let us
denote by $\Omega_{\alpha,\beta}(t,\lambda_1,\lambda_2)$ the value of
the symplectic pairing \eqref{phase-factor} at
$(t,\lambda_1,\lambda_2)\in D^+_\infty$. Using formulas
\eqref{phase-Wab} and \eqref{Wab-integral} we can extend analytically
and define $\Omega_{\alpha,\beta}(t,\lambda_1,\lambda_2)$ for all
$(t,\lambda_1)\in (H\times \CC)'$ and $\lambda_2$ sufficiently close
to $\lambda_1$ and $\lambda_2\neq \lambda_1$. 
Note however that the value of
$\Omega_{\alpha,\beta}(t,\lambda_1,\lambda_2)$, even if we fix a reference
path from $(t^\circ,\lambda^\circ)$ to $(t,\lambda_1)$, is determined only
up to an integer multiple of $2\pi\ii 
(\alpha|\beta)$ because, according to Proposition
\ref{ope-cal-periods}, the symplectic pairing has a logarithmic
singularity at $\lambda_1=\lambda_2$ of the form $(\alpha|\beta)\log
(\lambda_1-\lambda_2)$.
\begin{remark}
The definition of the domain $D^+_\infty$ in \cite{MST}, Section 4.4,
involves one more inequality $|\lambda_2|>r(t)$. This inequality is
necessary only if we take also the Laurent series expansion of
\eqref{phase-factor} at $\lambda_2=\infty$. 
\end{remark}

\subsection{Phase form}
The {\em phase 1-form} is by definition the 1-form 
\ben
\W_{\alpha,\beta}(t,\xi) = I^{(0)}_\alpha(t,0)\bullet I^{(0)}_\beta(t,\xi),
\een 
where the parameter $\xi$ is sufficiently small and the RHS is
interpreted as an element in $T^*_tH= H^*$ via the Poincare
pairing. Note that we are using the linear structure on $H$ to
trivialize the tangent bundle $TH=H\times H$. The form is defined for
all $t\in H'$, where $H'$ is defined to be the 
set of points where $E\bullet_t$ is invertible. Recall the
reference point $(t^\circ,\lambda^\circ)$ for $H\times \CC$ used in
the definition of the periods. We fix 
$t^\circ-\lambda^\circ \mathbf{1}\in H$ to be a reference point.   The value of the
form depends on the choice of a reference path in $H'$.   

Suppose that $C\subset H'$ is a simple loop based at the point
$t-\lambda_1\mathbf{1}\in H'$ that goes around a generic point on the
discriminant. Let $\varphi\in \mathfrak{R}$ be the reflection vector
corresponding to the simple loop corresponding to the composition of
$C$ and a reference path from $t^\circ-\lambda^\circ\mathbf{1}$ to
$t-\lambda_1\mathbf{1}$. Then the following formula holds:
\beq\label{phase-form-periods}
\Omega_{w(\alpha),w(\beta)}(t,\lambda_1,\lambda_2)-
\Omega_{\alpha,\beta}(t,\lambda_1,\lambda_2) =
2\pi\ii (\alpha|\varphi)\, \langle \varphi,\beta\rangle + 
\int_{t'\in C} \W_{\alpha,\beta}(t',\lambda_2-\lambda_1),
\eeq
where $w$ is the monodromy transformation corresponding to the closed
loop $C$ and $\lambda_2$ is sufficiently close to $\lambda_1$. Here
sufficiently close means the following: the line segment 
$[(t',\lambda_1-\lambda_2),(t',0)]$, when $t'$ varies along $C$, sweeps
out a fat loop in $H\times \CC$, which does not intersect the discriminant.
The proof of \eqref{phase-form-periods} is given in \cite{Mi2},
Section 4.4 in the settings of singularity theory. It generalizes in
our case too using results from \cite{Mi4} and \cite{MST}. For the
reader's convenience, we give a self-contained proof. 

The following result is proved in \cite{Mi4}, Lemma 15. 
\begin{proposition}\label{pf:euler-equation}
The phase form is weighted-homogeneous of weight $0$, i.e., 
\begin{eqnarray*}
(\xi\partial_\xi + L_E) \mathcal{W}_{\alpha,\beta}(s,\xi) = 0,
\end{eqnarray*} 
where $L_E$ is the Lie derivative with respect to the vector field $E$.
\end{proposition}
\proof
Note that 
\begin{eqnarray*}
\mathcal{W}_{\alpha,\beta}(s,\xi) = (d I^{(-1)}_\alpha(s,0), I^{(0)}_\beta(s,\xi)).
\end{eqnarray*}
It is easy to check that $\mathcal{W}_{\alpha,\beta} $ is a closed
1-form, so using the Cartan's magic formula $L_E=d\iota_E + \iota_Ed$,
where $\iota_E$ is the contraction by the vector field $E$, we get 
\begin{eqnarray*}
  L_E \mathcal{W}_{\alpha,\beta} =
  d( (\theta+1/2) I^{(-1)}_\alpha(s,0), I^{(0)}_\beta(s,\xi)) =
  -d( I^{(-1)}_\alpha(s,0), (\theta-1/2) I^{(0)}_\beta(s,\xi) ).
\end{eqnarray*}
We used that $\theta$ is skew-symmetric with respect to the residue pairing
and that 
\begin{eqnarray*}
\iota_E d I^{(-1)}_\alpha(s,0) = E I^{(-1)}_\alpha(s,0)) = (\theta+1/2) I^{(-1)}_\alpha(s,0),
\end{eqnarray*}
where the last equality is obtained by substituting $\lambda=0$ in the
differential equation $\nabla^{(-1)}_{\partial/\partial \lambda} I^{(-1)}_\alpha(t,\lambda)
= 0$. Furthermore, using the Leibnitz rule we get
\begin{eqnarray}\label{LieW}
  L_E \W_{\alpha,\beta}=
-(d I^{(-1)}_\alpha(s,0) , (\theta-1/2) I^{(0)}_\beta(s,\xi))
-( I^{(-1)}_\alpha(s,0)) , (\theta-1/2) d I^{(0)}_\beta(s,\xi). 
\end{eqnarray}
The first pairing on the RHS of \eqref{LieW}, using the differential equation
$(\theta-1/2) I^{(0)}(s,\xi) =
(\xi\partial_\xi+E) I^{(0)}(s,\xi) $, becomes
\begin{equation}\label{res-2}
(dI^{(-1)}_\alpha(s,0) , (\xi\partial_\xi + E)I^{(0)}_\beta(s,\xi)) =
\xi\partial_\xi\mathcal{W}_{\alpha,\beta}(s,\xi)
+(A I^{(0)}_\alpha(s,0)) , E\bullet I^{(1)}_\beta(s,\xi) ),
\end{equation}
where $A=\sum_i (\phi_i\bullet)dt_i$ and we used that the periods are
solutions to the second structure connection. 
Similarly, using the skew-symmetry of $\theta$ and the
differential equation $dI^{(0)}_\beta = - AI^{(1)}_\beta$, we
transform the 2nd pairing on the RHS of \eqref{LieW} into
\begin{equation}\label{res-1}
((\theta+1/2)  I^{(-1)}_\alpha(s,0) , A I^{(1)}_\beta(s,\xi))  = 
-(E\bullet I^{(0)}_\beta(s,0), A I^{(1)}_\alpha(s,\xi)).
\end{equation}
On the other hand, since
the Frobenius multiplication is commutative, the commutator 
$[A,E\bullet]=0$, so the terms (\ref{res-1}) and (\ref{res-2}) add up
to $ \xi\partial_\xi\mathcal{W}_{\alpha,\beta}(s,\xi).$ 
\qed

Proposition \ref{pf:euler-equation} yields the following identity
\beq\label{partialW}
\partial_{\lambda_2} \W_{\alpha,\beta}(t', \lambda_2-\lambda_1) =
d_{t'} \Big( \frac{1}{\lambda_2-\lambda_1} \left(
  (\theta+1/2) I^{(-1)}_\alpha(t',0),
  I^{(0)}_\beta(t',\lambda_2-\lambda_1) \right) \Big),
\eeq
where $d_{t'}$ is the de Rham differential on $H$.
\begin{proposition}\label{prop:partial-phase}
  The following formula holds:
  \ben
  \partial_{\lambda_1} \Omega_{\alpha,\beta}(t,\lambda_1,\lambda_2) =
  \frac{1}{\lambda_1-\lambda_2}
(I^{(0)}_\alpha(t,\lambda_1), (\lambda_2-E\bullet) I^{(0)}_\beta(t,\lambda_2)).
  \een
\end{proposition}
\proof
The de Rham differentials with respect to $t$ of both sides coincide
with the 1-form $\partial_{\lambda_1}
\W_{\alpha,\beta}(t-\lambda_1\mathbf{1},\lambda_2-\lambda_1)$, where
for the LHS we used formulas \eqref{phase-Wab} and
\eqref{Wab-integral} and for the RHS we used formula
\eqref{partialW}, the equation
$(\lambda_2-E\bullet)I^{(0)}_\beta(t,\lambda_2) = (\theta+1/2)
I^{(-1)}_\beta(t,\lambda_2)$, and the symmetry
$
\W_{\alpha,\beta}(t-\lambda_1\mathbf{1},\lambda_2-\lambda_1) =
\W_{\beta,\alpha}(t-\lambda_2\mathbf{1},\lambda_1-\lambda_2) 
. $
It is enough to prove the formula 
for $(t,\lambda_1,\lambda_2)\in D^+_\infty$, because the difference of
the values of both the RHS and the LHS at two different points
$(t^{(1)},\lambda_1^{(1)},\lambda_2^{(1)})$ and $(t^{(2)},\lambda_1^{(2)},\lambda_2^{(2)})$
with
$\lambda_2^{(1)}-\lambda_1^{(1)}=\lambda_2^{(2)}-\lambda_1^{(2)}=:\xi$
is given by the same path integral
$\int_{t^{(1)}-\lambda_1^{(1)}\mathbf{1}}^{t^{(2)}-\lambda_1^{(2)}\mathbf{1}}\partial_\xi
\W_{\alpha,\beta}(t',\xi)$.
Using formulas
\eqref{phase-Wab} and \eqref{Wab-integral}, we get
\ben
\lim_{t\to -\infty}  \partial_{\lambda_1}
\Omega_{\alpha,\beta}(t,\lambda_1,\lambda_2) =
\partial_{\lambda_1}
\Omega(
\widetilde{\ff}^+_\alpha(\lambda_1,z),
\widetilde{\ff}^-_\alpha(\lambda_2,z)) =
\frac{1}{\kappa} \sum_{s=1}^\kappa (\alpha|\sigma^s\beta)
\frac{\lambda_2^{1/\kappa} \lambda_1^{-1}}{ \eta^s
  \lambda_1^{1/\kappa}-\lambda_2^{1/\kappa} }.
\een
Note that
\ben
\lim_{t\to -\infty} \Big(
\frac{1}{\lambda_1-\lambda_2}
(I^{(0)}_\alpha(t,\lambda_1), (\lambda_2-E\bullet)
I^{(0)}_\beta(t,\lambda_2)) \Big) =
\frac{1}{\lambda_1-\lambda_2}
\left(
  \widetilde{I}^{\ (0)}_\alpha(\lambda_1), (\lambda_2-\rho)
  \widetilde{I}^{\ (0)}_\beta(\lambda_2)
\right).
\een
Therefore, we need only to prove that
\beq\label{phase-limit-id}
\frac{1}{\lambda_1-\lambda_2}
\left(
  \widetilde{I}^{\ (0)}_\alpha(\lambda_1), (\lambda_2-\rho)
  \widetilde{I}^{\ (0)}_\beta(\lambda_2)
\right) =
\frac{1}{\kappa} \sum_{s=1}^\kappa (\alpha|\sigma^s\beta)
\frac{\lambda_2^{1/\kappa} \lambda_1^{-1}}{ \eta^s
  \lambda_1^{1/\kappa}-\lambda_2^{1/\kappa} }.
\eeq
It is enough to check the above formula when $(\alpha,\beta)$ is a
pair of basis vectors, that is, $\alpha=\phi_{i,p}$ and
$\beta=\phi_{j,q}$.
To begin with, we need to find explicit formulas for the intersection
pairing and $\sigma$ in the basis $\{\phi_{i,p}\}$. 
The intersection pairing can be computed easily
using the formulas in Lemma \ref{le:ip}. Namely, we have
\ben
(\phi_{0,0}|\phi_{0,0}) = \frac{1}{n-2},\quad
(\phi_{i,p}|\phi_{i,a_i-p}) = \frac{\sin(\pi p/a_i)}{\pi a_i}\quad
(1\leq i\leq 3,\ 1\leq p\leq a_i-1)
\een
and all other pairings between the bases vectors are 0. In particular,
the kernel of the intersection pairing is spanned by
$\phi_{0,1}$. Since the classical monodromy operator $\sigma$ 
is defined by analytic continuation of
the periods in counterclockwise direction, using formula
\eqref{calibr-period}, we get 
\begin{align}
  \nonumber
\sigma(\phi_{0,0}) & = \phi_{0,0} + 2\pi\ii \, \frac{1}{n-2}\,
                    \phi_{0,1}, \\
  \nonumber
  \sigma(\phi_{0,1}) & =\phi_{0,1}, \\
  \nonumber
  \sigma(\phi_{i,p}) & = \eta_i^{-p} \phi_{i,p} \quad
  (1\leq i\leq 3,\ 1\leq p\leq a_i-1).
\end{align}
The proof of \eqref{phase-limit-id} is a straightforward
computation. Let us consider only the case when $\alpha=\phi_{i,p}$
and $\beta=\phi_{i,a_i-p}$ leaving the rest of the cases as an
exercise. The LHS of \eqref{phase-limit-id} is
\ben
\frac{1}{\lambda_1-\lambda_2}\, 
\frac{\lambda_1^{-p/a_i}}{\Gamma(1-p/a_i)} \,
\frac{\lambda_2^{p/a_i}}{\Gamma(p/a_i)}  =
\frac{\sin(\pi p/a_i)}{\pi a_i}
\frac{
(\lambda_2/\lambda_1)^{p/a_i}  }{\lambda_1-\lambda_2}.
\een
The RHS of \eqref{phase-limit-id}, since
$(\alpha|\sigma^s\beta) = \eta_i^{ps}(\phi_{i,p}|\phi_{i,a_i-p})$, is
\ben
\frac{\sin(\pi p/a_i)}{\pi a_i}\, \frac{1}{\kappa}
\sum_{s=1}^\kappa \eta_i^{ps}
\frac{\lambda_2^{1/\kappa} \lambda_1^{-1}}{ \eta^s
  \lambda_1^{1/\kappa}-\lambda_2^{1/\kappa} } =
\frac{\sin(\pi p/a_i)}{\pi a_i}\, \lambda_1^{-1}\, \frac{1}{\kappa}
\sum_{s=1}^\kappa
\sum_{m=1}^\infty
\eta^{(-m+p\kappa/a_i)s}(\lambda_2/\lambda_1)^{m/\kappa},
\een
where we took the Laurent series expansion at $\lambda_1=\infty$ and
used that $\eta_i^{ps} = \eta^{ps\kappa/a_i}$.  Let us exchange the
order of the sums over $s$ and $m$. The sum over $s$ is not zero if
and only if $m/\kappa = p/a_i +\ell$ for some integer $\ell\geq
0$. Therefore, the above sum, i.e., the RHS of \eqref{phase-limit-id}
is
\ben
\frac{\sin(\pi p/a_i)}{\pi a_i}\, \lambda_1^{-1}\,
\sum_{\ell=0}^\infty (\lambda_2/\lambda_1)^{\ell +p/a_i} =
\frac{\sin(\pi p/a_i)}{\pi a_i}\,
\frac{(\lambda_2/\lambda_1)^{p/a_i} }{ \lambda_1-\lambda_2},
\een
which is precisely the same as the LHS of \eqref{phase-limit-id}.
\qed

Under the same notation as in formula \eqref{phase-form-periods}, let
us introduce the following bilinear form:
\ben
I_C(\alpha,\beta):=
\int_{t'\in C} \W_{\alpha,\beta}(t',\lambda_2-\lambda_1) -
\Omega_{w(\alpha),w(\beta)}(t,\lambda_1,\lambda_2)+
\Omega_{\alpha,\beta}(t,\lambda_1,\lambda_2). 
\een
We have to prove that $I_C(\alpha,\beta) = - 2\pi\ii\,
(\alpha|\varphi)\, \langle\varphi,\beta\rangle$. 
Using the differential equations of the second
structure connection, it is easy to check that $I_C(\alpha,\beta)$ is
independent of $t$. Proposition \ref{prop:partial-phase} and formula
\eqref{partialW} imply that $I_C(\alpha,\beta)$ is also independent of
$\lambda_1$ and $\lambda_2$. Therefore, we may
assume that the loop $C$ is in a neighborhood of $t-u_i(t)\mathbf{1}$,
for some $1\leq i\leq n+1$, where $u_i(t)$ ($1\leq i\leq n+1$) is the set of canonical
coordinates, i.e., the eigenvalues of the operator $E\bullet_t$. For
homotopy reasons, we can also assume that the path $C$ has a
parametrization $t'=t+(x-\lambda_1)\mathbf{1}$, where the parameter
$x$ varies counterclockwise along a small loop $\gamma$ that encloses the points
$\lambda_1-u_i$ and $\lambda_2-u_i$. 

Suppose that formula 
\eqref{phase-form-periods} is proved in the special case when $\beta=\varphi$ and
$(\alpha|\varphi)=0$, that is,
\beq\label{paeriod-a_phi}
\int_{t'\in C} \W_{\alpha,\varphi}(t',\lambda_2-\lambda_1) =
-2\Omega_{\alpha,\varphi}(t,\lambda_1,\lambda_2).
\eeq
The general case can be deduced in the following way: 
Let us decompose
\ben
\alpha=:\alpha'+(\alpha|\varphi)\varphi/2,\quad \beta=:\beta'+(\beta|\varphi)\varphi/2,
\een
where $(\alpha'|\varphi)=(\beta'|\varphi)=0$. Note that
\ben
w(\alpha)=\alpha'-(\alpha|\varphi)\varphi/2,\quad
w(\beta)=:\beta'-(\beta|\varphi)\varphi/2.
\een
The quantity $I_C(\alpha,\beta)$ is bilinear in $\alpha$ and $\beta$, so
substituting the above formulas, after a short computation, we get
\begin{align}\label{index-1}
I_C(\alpha,\beta)= & \int_C \W_{\alpha',\beta'} (t',\lambda_2-\lambda_1) +
\tfrac{1}{4} (\alpha|\varphi)\, (\beta|\varphi) \, \int_{C}
  \W_{\varphi,\varphi}(t',\lambda_2-\lambda_1) +\\
  \label{index-2}
  &
\tfrac{1}{2}(\alpha|\varphi) \int_C
\W_{\varphi,\beta'}(t',\lambda_2-\lambda_1) +
\tfrac{1}{2}(\beta|\varphi) \int_C
    \W_{\alpha',\varphi}(t',\lambda_2-\lambda_1) + \\
  \label{index-3}
                 &
                   (\alpha|\varphi)
                   \Omega_{\varphi,\beta'}(t,\lambda_1,\lambda_2) +
                   (\beta|\varphi)
                   \Omega_{\alpha',\varphi}(t,\lambda_1,\lambda_2). 
\end{align}
The first integral in \eqref{index-1} is 0 because the cycles
$\alpha'$ and $\beta'$ are invariant with respect to the local
monodromy, which is equivalent to the fact that the periods
$I^{(0)}_{\alpha'}(t',0)$ and
$I^{(0)}_{\beta'}(t',\lambda_2-\lambda_1)$ are analytic in $t'$ in a
neighborhood of the point $t-u_i(t)\mathbf{1}$. 
The second integral in \eqref{index-1}, since it
coincides with $I_C(\varphi,\varphi)$, does not depend
on $(t,\lambda_1,\lambda_2)$, so we may set $\lambda_2=\lambda_1$. By definition, 
\ben
\int_C \W_{\varphi,\varphi}(t',0)  =
\int_{\gamma}
(I^{(0)}_\varphi(t,\lambda_1-x), I^{(0)}_\varphi(t,\lambda_1-x)),
\een
where $\gamma$ is the small loop in the $x$-plane
mentioned above based at $x=0$ and enclosing $\lambda_1-u_i$ and
$\lambda_2-u_i$. 
On the other hand, the period has the following expansion
\ben
I^{(0)}_\varphi(t',\xi) = \frac{\sqrt{2}}{\sqrt{\xi -u_i(t') }}
\Big(e_i + O(\xi-u_i)\Big),  
\een
where $e_i= \tfrac{du_i}{\sqrt{\Delta_i}}$ are the normalized
idempotent, that is, $(e_i,e_j)=\delta_{i,j}$. Clearly, only the
leading order term contributes, so
\beq\label{phase-period}
\int_C \W_{\varphi,\varphi}(t',0)  = \int_{\gamma}
\frac{2dx}{\lambda_1-u_i-x} = -4\pi\ii. 
\eeq
Both integrals in \eqref{index-2} can be computed with formula
\eqref{paeriod-a_phi}, that is, the first integral is
$-2\Omega_{\beta',\varphi}(t,\lambda_2,\lambda_1)$ and the second one
is $-2\Omega_{\alpha',\varphi}(t,\lambda_1,\lambda_2)$. We get the
following formula:
\ben
I_C(\alpha,\beta) = -\pi\ii (\alpha|\varphi)(\beta|\varphi) +
(\alpha |\varphi) \Big(
\Omega_{\varphi,\beta'}(t,\lambda_1,\lambda_2) -
\Omega_{\beta',\varphi}(t,\lambda_2,\lambda_1)
\Big).
\een
Recalling formulas \eqref{phase-Wab} and \eqref{calibr-locality}, we
finally get
\ben
I_C(\alpha,\beta) = -\pi\ii (\alpha|\varphi)(\beta|\varphi)-2\pi\ii
(\alpha|\varphi)\langle \varphi,\beta'\rangle. 
\een
Substituting $\beta'=\beta-(\beta|\varphi) \, \varphi/2$ and using
$\langle \varphi,\varphi\rangle =1$, we get $I_C(\alpha,\beta)= -2\pi\ii
(\alpha|\varphi)\langle \varphi,\beta \rangle$, which is what we had
to prove.

For the proof of \eqref{paeriod-a_phi}, the most delicate step is to
reduce the proof to the case when the loop $C$ belongs to a
neighborhood of a point $t-u_i(t)\mathbf{1}$, such that
$|u_i(t)|>|u_j(t)|$ for all $i\neq j$. In order to do this, we need to
use the {\em homotopy invariance property} of $I_C(\alpha,\beta)$ and
the so-called {\em Painleve property} of a semi-simple Frobenius
manifold (see \cite{Du2}). Let us state the homotopy invariance property. Suppose that
$h:[0,1]\times [0,1]\to H'$ is a continuous map, such that
$\gamma_s=h(s,\ ):[0,1]\to H'$ is a simple loop around the discriminant
for all $s\in [0,1]$. Let $C_0=\gamma_0$ and $C_1=B^{-1}\circ \gamma_1\circ
B$, where $B$ is the path $B(s) = h(0,s)$ for $s\in [0,1]$. Then
$I_{C_0}(\alpha,\beta) = I_{C_1}(\alpha,\beta)$. For the proof, we
need only to notice that the difference
\ben
I_{C_0}(\alpha,\beta) - I_{C_1}(\alpha,\beta) = \int_{[0,1]\times
  [0,1]} h^*(d\W_{\alpha,\beta}).
\een
It is easy to check that the phase form $\W_{\alpha,\beta}$ is closed,
so the above integral must be $0$. Let us state the Painleve property.
Let us assume that $t$ is a 
semi-simple point, such that  the canonical coordinates give an embedding of a
neighborhood $U$ of $t$ in the configuration space
\ben
\mathfrak{M}_{n+1}:= \{u\in \CC^{n+1}\ |\ u_i\neq u_j\ \mbox{ for } \ i\neq j \}.
\een
The Painleve property says that the Frobenius structure on $U\subset
\mathfrak{M}_{n+1}$ can be extended along any path in $\mathfrak{M}_{n+1}$ except for a set
of paths that correspond to the points of an analytic hypersurface in
the universal cover of $\mathfrak{M}_{n+1}$. The Painleve
property is very difficult to prove. It is a corollary of a result of Malgrange about the analytic
continuation of the solutions of the Schlesinger equations (see \cite{Mal}). We refer
to  \cite{Mi3} for a self-contained proof and for further
references. Therefore, for every given $i$ ($1\leq i\leq n+1$), we can 
choose a generic point $u'\in \mathfrak{M}_{n+1}$, such that the $i$-th
coordinate $u'_i$ of $u'$ has 
the largest absolute value and there exists a path $A:[0,1]\to
\mathfrak{M}_{n+1}$ connecting $u(t)$ and
$u'$, such that the Frobenius structure extends along it. Clearly, if
we have a simple loop $\gamma\subset \CC$, based at $\lambda_1$, around 
$u_i(t)=A_i(0)$, we can extend it to a continuous family $\gamma_s$
($0\leq s\leq 1$), such that $\gamma_s\subset \CC$ is a simple loop
around $A_i(s)$. Therefore, $h(s_1,s_2):=
A(s_1)-\gamma_{s_1}(s_2)\mathbf{1}$ is a homotopy which allows us to
replace the simple loop $C(s):=t-\gamma(s)\mathbf{1}$ around
$t-u_i(t)\mathbf{1}$ with a simple loop around the point
$t'-u_i(t')\mathbf{1}$, where $t'=A(1)$. By construction $u_i(t') =
A_i(1) = u'_i$ has the largest absolute value among all coordinates of
$u'$. 

Let us prove \eqref{paeriod-a_phi} assuming the following conditions:  
the path $C$ is a sufficiently small loop in
a neighborhood of a generic point on the discriminant with
parametrization $t'=t+(x-\lambda_1)\mathbf{1}$, where $x$ varies
along a small loop $\gamma$ based at $0$ enclosing $\lambda_1-u_i$ and 
$\lambda_2-u_i$, and $|u_i(t)|>|u_j(t)|$ for all $i\neq j$.
We still have the freedom to choose $\lambda_1$ and $\lambda_2$ as we
wish provided that they are sufficiently close to $u_i(t)$. First, we
choose a small disk $\Delta$ with center $0$, such that,
\ben
|u_i(t)|>|u_j(t)+x|,\quad \forall x\in \Delta,\quad
\forall j\neq i.
\een
Let us choose $\lambda_1$ and $\lambda_2$ to be on the line passing
through $0$ and $u_i(t)$, such that, 
$|\lambda_1|>|\lambda_2|>|u_i(t)|$ and $\lambda_1-u_i(t), \lambda_2-u_i(t)\in
\Delta$. Note that $(t,\lambda_1,\lambda_2)\in D^+_\infty$ and that
for homotopy reasons we may assume that the loop $\gamma\subset
\Delta$. Using integration by parts and that
$I^{(-n-1)}_\varphi(t,\lambda_2-x)$ vanishes at  
$x=\lambda_2-u_i$, we get for all $x\in \Delta$
\beq\label{phase-ui}
\int_{\lambda_2-u_i}^x (I^{(0)}_\alpha(t,\lambda_1-y),
I^{(0)}_\varphi(t,\lambda_2-y)) dy =
\sum_{n=0}^\infty (-1)^{n+1}
(I^{(n)}_\alpha(t,\lambda_1-x),I^{(-n-1)}_\varphi(t,\lambda_2-x)),
\eeq
where the RHS is interpreted via its Laurent series expansion at
$\lambda_1=\infty$ and the integration is along the straight segment
$[\lambda_2-u_i,x]$. The period
$I^{(0)}_\alpha(t,\lambda_1-y)=I^{(0)}_\alpha(t+y\mathbf{1},\lambda_1)$ is
analytic at $\lambda_1=u_i(t)+y$, because the condition
$(\alpha|\varphi)=0$ implies that it is invariant with respect to
the local monodromy. Since all other singularities are at the points
$\lambda_1=u_j(t)+y$ for $j\neq i$ and $|u_j(t)+y|<|u_i(t)|<|\lambda_1|$, we get
that the Laurent series expansion at $\lambda_1=\infty$ of
$I^{(0)}_\alpha(t,\lambda_1-y) $ is convergent for all
$|\lambda_1|>|u_i(t)|$ and it depends analytically on $y\in
\Delta$. Since $I^{(0)}_\varphi(t,\lambda_2-y)$ has at most a pole of
order $\tfrac{1}{2}$ at $y=\lambda_2-u_i$, we get that the LHS of
\eqref{phase-ui} has the form  $(\lambda_2-u_i-x)^{1/2} f(x,\lambda_1)$ for some function $f$ holomorphic for
all $x\in \Delta$ and $|\lambda_1|>|u_i|$. Therefore, since the loop
$\gamma\subset \Delta$, we get
\ben
-2 (\lambda_2-u_i)^{1/2} f(0) =
\int_\gamma (I^{(0)}_\alpha(t,\lambda_1-x),
I^{(0)}_\varphi(t,\lambda_2-x)) dx =
\int_\gamma\W_{\alpha,\varphi}(t',\lambda_2-\lambda_1).
\een
It remains only to note that the RHS of \eqref{phase-ui} at $x=0$, since the point
$(t,\lambda_1,\lambda_2)\in D^+_\infty$, is precisely
$\Omega_{\alpha,\varphi}(t,\lambda_1,\lambda_2)$.

\section{The symplectic space formalism}

We would like to recall Givental's symplectic space formalism, which
we will need in the proof of Theorem \ref{t1}.
\subsection{Quantizing symplectic transformation}
Suppose that $A:\H\to \H$ is an infinitesimal symplectic
transformation
\ben
\Omega(Af,g)+\Omega(f,Ag)=0,\quad \forall f,g\in \H.
\een
If $\ff\in \H$, then let us decompose it as
\ben
\ff=\sum_{k=0}^\infty
\sum_{i,a} (
p_{i,a,k}\phi^{i,a}(-z)^{-k-1}+
q_{i,a,k}\phi_{i,a}z ^k),
\een
where $\{\phi_{i,a}\}$ is the basis of $H$ that we fixed already and
$\{\phi^{i,a}\}$ is its dual basis with respect to the Poincare
pairing. The numbers $p_{i,a,k}$ and $q_{i,a,k}$ can be viewed as
coordinates of $\ff$.
Under the natural trivialization of the tangent
bundle $T\H=\H\times \H$ the symplectic pairing $\Omega$ defines a
symplectic form on $\H$ which coincides with $\sum_{i,a,k}
dp_{i,a,k}\wedge dq_{i,a,k}$. In other words, $\H$ is a symplectic
manifold and $p_{i,a,k}$ and $q_{i,a,k}$ form a Darboux coordinate
system. The linear map $\ff\mapsto A\ff$ defines a linear
Hamiltonian vector field with Hamiltonian
\ben
h_A(\ff):=\frac{1}{2}\Omega(A\ff,\ff).
\een
We define $\widehat{A}:=\widehat{h}_A$ where quadratic functions are
quantized according to the following rules:
\ben
(q_{i,a,k}q_{j,b,\ell})\sphat :=\hbar^{-1}q_{i,a,k}q_{j,b,\ell},\quad
(p_{i,a,k}p_{j,b,\ell})\sphat = \hbar
\frac{\partial^2}{\partial q_{i,a,k}\partial q_{j,b,\ell}}.
\een
and
\ben
(q_{i,a,k} p_{j,b,\ell})\sphat =
(p_{j,b,\ell}q_{i,a,k})\sphat =
q_{i,a,k}\frac{\partial}{\partial q_{j,b,\ell}}, \quad
\een
Finally, if $M=e^A$ is a symplectic transformation then we define
$\widehat{M}:=e^{\widehat{A}}$.

\subsection{Givental's higher genus reconstruction}
Suppose that $t\in H$ is a semi-simple point. The eigenvalues
$\{u_i(t)\}_{1\leq i\leq N}$
of $E\bullet_t$ form a local coordinate system in which the Frobenius
pairing and multiplication take the form
\ben
({\partial}/{\partial u_i},{\partial}/{\partial u_j} )=
\delta_{i,j}/\Delta_j(t),\quad 
\frac{\partial}{\partial u_i}\bullet \frac{\partial}{\partial u_j}=
\delta_{i,j}
\frac{\partial}{\partial u_i}
\een
where $\Delta_j$ $(1\leq j\leq N)$ are functions such that
$\Delta_j(t)\neq 0$.
Let $e_i(t):=\sqrt{\Delta_i}\partial/\partial u_i\in T_tH=H$ be the
normalized idempotents. There exists a unique formal power series
\ben
R(t,z) = 1+\sum_{k=1}^\infty R_k(t) z^k,\quad
R_k(t)\in \operatorname{End}(H)
\een
such that the functions $R(t,z)e_i(t) e^{u_i(t)/z}$ are solutions to
the Dubrovin's connection (see \cite{G1}).

Let $\D_{\rm pt}(\hbar,\tt)$ be the total descendant potential of a
point (before the dilaton shift). Let us denote by $\D^{(i)}_{\rm
  pt}(\hbar,\qq)$ the formal series obtained from $\D_{\rm
  pt}(\hbar,\tt)$ via the substitutions
\ben
\hbar\mapsto \hbar \Delta_i,\quad
t_k:= q_k(u_i) + \delta_{k,1}
\een
where $q_k=\sum_{i,a}q_{i,a,k}\partial_{i,a}$ is identified with a
flat vector field acting on the function $u_i$ by derivation. 
According to Givental's higher genus reconstruction \cite{G1}, proved
by Teleman \cite{Te}, the total descendant potential of the orbifold
$\PP^1_{n-2,2,2}$ is given by the following formula
\ben
\D(\hbar,\qq) =
e^{F^{(1)}(t)} (S(t,Q,z)^{-1})\sphat\ 
(R(t,z))\sphat \ \prod_{i=1}^N \D^{(i)}_{\rm pt}(\hbar,\qq),
\een
where $F^{(1)}(t)$ is the genus-1 potential without descendants. It is
known that both operator series $S(t,Q,z)$ and 
$R(t,z)$ are symplectic, so their quantization makes sense. 

\subsection{Conjugation by $S$}
We already discussed the conjugation by $S$ of the vertex operators
(see formula \eqref{S-conj}).  We get 
\ben
(\widetilde{\Gamma}^\alpha(\lambda) \otimes
\widetilde{\Gamma}^{-\alpha} (\lambda))
(\widehat{S}^{-1}\otimes \widehat{S}^{-1}) =
(\widehat{S}^{-1}\otimes \widehat{S}^{-1})
e^{W_{\alpha,\alpha}(t,\lambda,\lambda)}
(\Gamma^\alpha(t,\lambda)\otimes
\Gamma^{-\alpha}(t,\lambda)).
\een
Using formula \eqref{phase-Wab}, we get
\beq\label{W-phase}
W_{\alpha,\beta}(t,\lambda,\lambda) =
\lim_{\mu\to \lambda} \Big(
\Omega_{\alpha,\beta}(t,\lambda,\mu)-
\widetilde{\Omega}_{\alpha,\beta}(\lambda,\mu)
\Big),\quad (t,\lambda,\lambda)\in D_\infty.
\eeq
Using the definition of $W_{k,\ell}$ we get
\ben
\partial_\lambda W_{\alpha,\alpha}(t,\lambda,\lambda) =
-(I^{(0)}_\alpha(t,\lambda),I^{(0)}_\alpha(t,\lambda))+
(\widetilde{I}^{(0)}_\alpha(\lambda),
\widetilde{I}^{(0)}_\alpha(\lambda)).
\een
Recalling the explicit formulas for the calibrated periods, Lemma
\ref{le:a-euler-b0}, and Lemma \ref{le:ab_tw} we get
\ben
\partial_\lambda W_{\alpha,\alpha}(t,\lambda,\lambda) =
-(I^{(0)}_\alpha(t,\lambda),I^{(0)}_\alpha(t,\lambda))+
\Big((\alpha|\alpha) +\frac{1}{n-2}
\operatorname{rk}(\alpha)
\operatorname{rk}(\alpha)\Big)\lambda^{-1} .
\een

\subsection{Conjugation by $R$}\label{sec:conj-R}
Let $(t,\lambda)$ be a
point in a neighborhood of a generic point on the discriminant. In
other words $\lambda$ is sufficiently close to $u:=u_i(t)$. Let
$\beta=\epsilon'-\epsilon'' $ ($\epsilon',\epsilon''\in \E$ ) be the
corresponding vanishing cycle. Put
$
\alpha:=(\epsilon'+\epsilon'')/2.
$
Then we have the following factorizations (see \cite{G3}, Proposition 4)
\ben
\Gamma^{\epsilon'}(t,\lambda)\otimes\Gamma^{-\epsilon'}(t,\lambda) = 
e^{K_{\alpha,\beta}(t,\lambda)} (
\Gamma^{\alpha}(t,\lambda)\otimes\Gamma^{-\alpha}(t,\lambda))
(\Gamma^{\beta/2}(t,\lambda)\otimes\Gamma^{-\beta/2}(t,\lambda)  ),
\een
and
\ben
\Gamma^{\epsilon''}(t,\lambda)\otimes\Gamma^{-\epsilon''}(t,\lambda) = 
e^{-K_{\alpha,\beta}(t,\lambda)} (
\Gamma^{\alpha}(t,\lambda)\otimes\Gamma^{-\alpha}(t,\lambda))
(\Gamma^{-\beta/2}(t,\lambda)\otimes\Gamma^{\beta/2}(t,\lambda)  ),
\een
where
\ben
K_{\alpha,\beta}(t,\lambda) =
-\int_{t-u\mathbf{1}}^{t-\lambda\mathbf{1}} \W_{\alpha,\beta},
\een
where the integration path is along the straight segment between the
points  $t-\lambda\mathbf{1}$ and $t-u\mathbf{1}$.

Let us define the vertex operators
\ben
\Gamma^{\pm }_i:=
\exp\Big(\sum_{n\in \ZZ} (-z\partial_\lambda)^n
\frac{e_i}{\sqrt{2(\lambda-u_i)} }\Big)\sphat,
\een
where $e_i=\sqrt{\Delta_i}\partial/\partial u_i$ are the normalized
idempotents.   Then the following formula holds (see \cite{G3}, page 490):
\ben
(\Gamma^{\pm \beta/2}\otimes \Gamma^{\mp \beta/2})
(\widehat{R}\otimes \widehat{R}) =
e^{V_{\beta,\beta}}
(\widehat{R}\otimes \widehat{R})
(\Gamma^{\pm}_i\otimes \Gamma^{\mp}_i),
\een
where
\ben
V_{\beta,\beta} =- \lim_{\epsilon \to 0}
\int_{t-(u_i+\epsilon)\mathbf{1}}^{t-\lambda\mathbf{1}}
\Big(
\W_{\beta/2,\beta/2} + \frac{du_i}{2u_i}\Big),
\een
where the limit is taken along a straight segment connecting
$t-\lambda\mathbf{1}$ and $t-u_i\mathbf{1}$. 

Using the factorization and the conjugation by $R$ formulas we get the
following important identity
\ben
(\Gamma^\varphi\otimes
\Gamma^{-\varphi})(\widehat{R}\otimes\widehat{R}) =
\frac{c^i_\alpha(t,\lambda)}{c^i_\varphi(t,\lambda)}
(\Gamma^{{\rm sign}(\alpha|\varphi)\alpha}\otimes
\Gamma^{-{\rm sign}(\alpha|\varphi)\alpha})
(\widehat{R}\otimes\widehat{R})
(\Gamma_i^{{\rm sign}(\beta|\varphi)}\otimes
\Gamma_i^{-{\rm sign}(\beta|\varphi)}),
\een
where ${\rm sign}(x)$ is the sign of the number $x$, $\varphi\in \{\pm
\epsilon',\pm\epsilon''\}$, and the coefficients $c^i_a(t,\lambda)$
for $a\in K^0(\PP^1_{n-2,2,2})$ are defined by
\ben
c_a^i(t,\lambda):=\lim_{\epsilon\to 0} \exp\Big(
\int_{t-(u_i+\epsilon)\mathbf{1}}^{t-\lambda\mathbf{1}}
\Big(
\W_{a,a} + (a|\beta)^2 \frac{du}{2u}
\Big)\Big).
\een
\subsection{Proof of Theorem \ref{t1}}

The proof is similar to the proof of Theorem 1.5 in \cite{MT}. Let us
write the total descendant potential as
$\D(\hbar,\qq)=e^{F^{(1)}(t)}\widehat{S(t,Q,z)}^{-1}\A_t(\hbar;\qq)$. The
formal series $\A_t$ coincides with the so called {\em total ancestor
  potential} of $\PP^1_{n-2,2,2}$ (see \cite{G1,G2}). Substituting
this formula in the Hirota quadratic equations of the extended D-Toda
hierarchy and conjugating with $\widehat{S}^{-1}$ we reduce the proof
to proving that the total ancestor potential satisfies the following
Hirota quadratic equations: for every $m\in \ZZ$ the 1-form
\beq\label{hqe-t}
\frac{d\lambda}{\lambda}
\Big(\Gamma^\#(t,\lambda)\otimes \Gamma(t,\lambda)\Big)
\Big(\sum_{\epsilon\in \E} 
b_{\epsilon}(t,\lambda) 
\Gamma^\epsilon(t,\lambda)\otimes
\Gamma^{-\epsilon}(t,\lambda)
\Big)(\tau\otimes \tau)
\eeq
computed at
\beq\label{discr-t}
\Omega(w_t,\qq'-\qq'')=m\hbar^{1/2}
\eeq
is regular in $\lambda$. Here the notation is as follows
\ben
w_t:=S(t,Q,z)\phi^{0,0}(-z)^{-1},
\een
the vertex operator (with coefficients in the algebra of differential operators)
\ben
\Gamma(t,\lambda) :=
\exp\Big( (
\ff_\varphi(t,\lambda)-w_t)\hbar^{1/2}\partial_x\Big)\sphat\ 
\exp\Big(\hbar^{-1/2}x v_t\Big)\sphat, 
\een
where $v_t:=S(t,Q,z) \phi_{0,0}$ and $\varphi:=1-L$, and the
coefficients
\ben
b_\epsilon(t,\lambda)=\widetilde{b}_\epsilon(\lambda)
e^{W_{\epsilon,\epsilon}(t,\lambda,\lambda)},
\een
where the notation $W_{\epsilon,\epsilon}$ is the same as in
\eqref{Wab}. We claim that the coefficients $b_\epsilon(t,\lambda)$
are compatible with the monodromy representation, i.e., the analytic
continuation of $b_{\epsilon}(t,\lambda)$ along a closed loop $C$
around the discriminant coincides with $b_{w(\epsilon)}(t,\lambda)$,
where $w$ is the monodromy transformation corresponding to the loop
$C$. Let us first prove the following simple lemma.
\begin{lemma}\label{le:calibr-b}
The coefficient $\widetilde{b}_\epsilon(\lambda)$ can be computed by
the following formula:
\ben
\frac{\lambda}{ \widetilde{b}_\epsilon(\lambda) } =
\lim_{\mu\to \lambda}\ \Big(
(\mu-\lambda)\, 
e^{\widetilde{\Omega}_{\epsilon,-\epsilon}(\lambda,\mu)}
\Big)\
e^{2\pi\ii \langle\epsilon,\epsilon^3_1\rangle }.
\een
\end{lemma}
\proof
Let us first point out that, thanks to the explicit
formulas in Section \ref{sec:eup}, the factor
\ben
e^{2\pi\ii \langle\epsilon,\epsilon^3_1\rangle } =
\begin{cases}
  -1 & \mbox{ if } \epsilon=\pm \epsilon^3_1, \\
  1 & \mbox{ otherwise.}
\end{cases}
\een
The rest of the proof is an explicit computation using Proposition
\ref{ope-cal-periods}. Let us consider first the case
$\epsilon=\epsilon^3_1$. We have
$\operatorname{rk}(\epsilon)=\operatorname{deg}(\epsilon)=0$ and
$\sigma(\epsilon)=-\epsilon$. Therefore, according to Proposition
\ref{ope-cal-periods}, we have
\ben
e^{\widetilde{\Omega}_{\epsilon,-\epsilon}(\lambda,\mu)} =
\prod_{s=1}^{\kappa} (1-\eta^{-s}
(\mu/\lambda)^{1/\kappa})^{(-1)^{s+1} } =
\frac{1-\tfrac{\mu}{\lambda} }{
 \left( 1-\tfrac{\sqrt{\mu}}{\sqrt{\lambda} } \right)^2} = 
\frac{\lambda-\mu}{(\sqrt{\lambda}-\sqrt{\mu})^2},
\een
where we used that $\kappa=2(n-2)$. Using the above formula we get
that the limit in the identity that we want to prove is $4\lambda =
\lambda/\widetilde{b}_{\epsilon^3_1}(\lambda)$. The computation in the
case $\epsilon=\epsilon^2_1$ is identical. Finally, suppose that
$\epsilon=\epsilon^1_i$. By definition
$\operatorname{rk}(\epsilon)=1$ and 
$\operatorname{deg}(\epsilon)=\tfrac{1}{2}+\tfrac{i}{n-2}$. Recalling Proposition
\ref{pr:cal-periods} we get that the pairing
$(\epsilon^1_i|\sigma^s(\epsilon^1_i))=1$ only if $s$ is an integer
multiple of $n-2$, and it is otherwise equal to $0$. Therefore,
according to Proposition \ref{ope-cal-periods}, we have
\ben
e^{\widetilde{\Omega}_{\epsilon,-\epsilon}(\lambda,\mu)} =
-\eta_1^i \, \frac{\mu^{1/(n-2)}}{Q}\,
\left(1-(\mu/\lambda)^{1/(n-2)}\right)^{-1} = -
\frac{\eta_1^i }
{Q( \mu^{-1/(n-2)}-\lambda^{-1/(n-2)} )} 
.
\een
Using the above formula we compute
\ben
\lim_{\mu\to \lambda}\ \Big( (\mu-\lambda)
e^{\widetilde{\Omega}_{\epsilon,-\epsilon}(\lambda,\mu)} \Big)=
\frac{(n-2) \eta_1^i }{Q}\, \lambda^{1+1/(n-2)} =
\frac{\lambda}{\widetilde{b}_{\epsilon^1_i}(\lambda)}.
\qed
\een
Recalling \eqref{phase-Wab} we get the following formula:
\beq\label{ope-b}
\frac{\lambda}{ b_\epsilon(t,\lambda) } =
\lim_{\mu\to \lambda}\ \Big(
(\mu-\lambda)\, 
e^{\Omega_{\epsilon,-\epsilon}(t,\lambda,\mu)}
\Big)\
e^{2\pi\ii \langle\epsilon,\epsilon^3_1\rangle }.
\eeq
Now we can prove the claim about the monodromy invariance. Suppose
that $C$ is a simple loop around the discriminant and that $\alpha$ is the
corresponding reflection vector. The analytic continuation of
$\lambda/b_{\epsilon}(t,\lambda)$ along the path $C$ is given by
\ben
\lim_{\mu\to \lambda}\ 
(\mu-\lambda)\, 
\exp\left( \Omega_{\epsilon,-\epsilon}(t,\lambda,\mu)+
  \int_C \W_{\epsilon,-\epsilon}(t',\mu-\lambda)\right)
\
e^{2\pi\ii \langle\epsilon,\epsilon^3_1\rangle }.
\een
We have to prove that the above limit coincides with
\ben
\lim_{\mu\to \lambda}\ 
(\mu-\lambda)\, 
\exp\left( \Omega_{w(\epsilon),-w(\epsilon)}(t,\lambda,\mu)\right)
\
e^{2\pi\ii \langle w(\epsilon),\epsilon^3_1\rangle },
\een
where $w(x)=x-(\alpha|x)\alpha$ is the reflection representing the
monodromy transformation along the loop $C$. In other words we have to
prove that
\beq\label{mon-inv-id}
\exp\Big(
I_C(\epsilon,-\epsilon) +
2\pi\ii \langle\epsilon,\epsilon^3_1\rangle-
2\pi\ii \langle w(\epsilon),\epsilon^3_1\rangle\Big)=1.
\eeq
According to formula \eqref{phase-form-periods},
$I_C(\epsilon-\epsilon)=2\pi\ii (\alpha|\epsilon)\langle
\alpha,\epsilon\rangle$, so 
\ben
I_C(\epsilon,-\epsilon) +
2\pi\ii \langle\epsilon,\epsilon^3_1\rangle-
2\pi\ii \langle w(\epsilon),\epsilon^3_1\rangle =
2\pi\ii\, (\alpha|\epsilon)\, \langle\alpha,\epsilon^3_1+\epsilon\rangle.
\een
Both $(\alpha|\epsilon)$ and
$\langle\alpha,\epsilon^3_1+\epsilon\rangle$ are integers as one can
see directly from our explicit formulas for the Euler pairing in
Section \ref{sec:eup}, so the identity in \eqref{mon-inv-id} is true
and the claim about the monodromy invariance of the coefficients
$b_\epsilon(t,\lambda)$ is complete.

Note that the 1-form \eqref{hqe-t} is invariant under the entire
monodromy group. Indeed the monodromy transformations up to
translations by terms of the type  $r(1-L)$ ($r\in \ZZ$) act on the set
$\E$ via permutations. Since
\ben
(
\Gamma(t,\lambda)^\#\otimes \Gamma(t,\lambda))
(\Gamma^{r\varphi}\otimes \Gamma^{-r\varphi}) =
e^{2\pi\ii\, \Omega(w_t,\qq'-\qq'')\, \hbar^{-1/2}}
\Gamma(t,\lambda)^\#\otimes \Gamma(t,\lambda)
\een
the substitution \eqref{discr-t}  eliminates the contributions of the
translation terms. The monodromy invariance implies that the 1-form
\eqref{hqe-t} is a formal power series in $\qq'+\mathbf{1}$ and
$\qq''+\mathbf{1}$ whose coefficients are formal Laurent series in
$\hbar^{1/2}$, whose coefficients are formal Laurent series in
$\lambda^{-1}$. By definition the regularity condition means
that the Laurent series must have only positive powers of $\lambda$,
i.e., it is polynomial. On the other hand the total ancestor potential
is a {\em tame} asymptotical function (see \cite{G3}, Proposition
5). This implies that the coefficient in front of each monomial in
$\hbar^{1/2}$, $\qq'+\mathbf{1}$, and $\qq''+\mathbf{1}$ must be a
polynomial expression in the period vectors
$I^{(m)}_{\epsilon}(t,\lambda)$ ($m\in \ZZ$, $\epsilon\in
\E$). Therefore each coefficient is a meromorphic function in
$\lambda\in \CC$ with possible poles only at the canonical coordinates
$u_i(t)$. In order to prove the polynomiality we need to check that
the 1-form \eqref{hqe-t} does not have a pole at $\lambda=u_i$. 

Let us pick $\lambda$ sufficiently close to $u:=u_i(t)$ and fix a
reference path that determines the values of all period
vectors. Suppose that $\beta:=\epsilon'-\epsilon''$ is a reflection
vector corresponding to the simple loop  around $\lambda=u_i$. The
only terms that could have a pole at $\lambda=u_i$ are
\beq\label{term-1}
\frac{d\lambda}{\lambda}
\Big(b_{\epsilon'}(t,\lambda)
\Gamma^{\epsilon'}(t,\lambda)\otimes \Gamma^{-\epsilon'}(t,\lambda)
+
b_{\epsilon''}(t,\lambda)
\Gamma^{\epsilon''}(t,\lambda)\otimes \Gamma^{-\epsilon''}(t,\lambda)
\Big)\A_t\otimes \A_t
\eeq
and
\ben
\frac{d\lambda}{\lambda}
\Big(
b_{\epsilon'}(t,\lambda)
\Gamma^{-\epsilon'}(t,\lambda)\otimes \Gamma^{\epsilon'}(t,\lambda)
+
b_{\epsilon''}(t,\lambda)
\Gamma^{-\epsilon''}(t,\lambda)\otimes \Gamma^{\epsilon''}(t,\lambda)
\Big)
\A_t\otimes \A_t.
\een
We claim that both are regular at $\lambda=u_i$. Let us prove this
statement for the first term. The argument for the second term is
similar. Let us compose the vertex operator expression in
\eqref{term-1} with $\widehat{R}\otimes 
\widehat{R}$. Recalling the formulas from Section \ref{sec:conj-R} we
get the composition of 
\beq\label{reg-term}
c^i_\alpha(t,\lambda) \Gamma^\alpha(t,\lambda)\otimes
\Gamma^{-\alpha}(t,\lambda)
(\widehat{R}\otimes \widehat{R})
\eeq
and
\beq\label{kdv-term}
\frac{d\lambda}{\lambda}
\Big(
\frac{b_{\epsilon'}(t,\lambda)}{c^i_{\epsilon'}(t,\lambda)}
\Gamma^+_i\otimes \Gamma^-_i +
\frac{b_{\epsilon''}(t,\lambda)}{c^i_{\epsilon''}(t,\lambda)}
\Gamma^+_i\otimes \Gamma^-_i
\Big)
\eeq
Note that the term \eqref{reg-term} is regular at $\lambda=u_i$,
because $\alpha$ is invariant with respect to the local monodromy
around $\lambda=u_i$, so all periods $I^{(m)}_\alpha(t,\lambda)$ are
analytic at $\lambda=u_i$. 
In order to finish the proof we just need to check that up to a
constant the operator \eqref{kdv-term} coincides with
\ben
\frac{d\lambda}{\sqrt{\lambda-u_i}}
\Big(
\Gamma^+_i\otimes \Gamma^-_i -
\Gamma^+_i\otimes \Gamma^-_i
\Big).
\een
If we prove this fact then according to Givental \cite{G3} the above
operator is the Hirota bilinear operator that defines the Hirota
bilinear equations for the KdV hierarchy, i.e., if we apply it to
$\tau\otimes\tau$ where $\tau$ is a tau-function of the KdV then we
get an expression regular in $\lambda$. The regularity that we would
like to prove follows because $\A_t=\widehat{R}\prod_i
\D^{(i)}_{\rm pt}$ and $\D^{(i)}_{\rm pt}$ is a tau-function of KdV as
it was conjectured by Witten \cite{Wi} and proved by Kontsevich \cite{Ko}.

Let us first prove that
\ben
\frac{b_{\epsilon'}(t,\lambda)}{c^i_{\epsilon'}(t,\lambda)} = A'
\frac{\lambda}{\sqrt{\lambda -u_i}}
\een
and
\ben
\frac{b_{\epsilon''}(t,\lambda)}{c^i_{\epsilon''}(t,\lambda)} = A''
\frac{\lambda}{\sqrt{\lambda -u_i}}
\een
for some constants $A'$ and $A''$. Note that
\ben
\partial_\lambda W_{\epsilon,\epsilon}(t,\lambda,\lambda) =
-(I^{(0)}_\epsilon(t,\lambda), I^{(0)}_\epsilon(t,\lambda)) +
\Big(1+\frac{1}{n-2} \operatorname{rk}(\epsilon)^2\Big)\lambda^{-1}
\een
and
\ben
\partial_\lambda \log c^i_\epsilon = -
(I^{(0)}_\epsilon(t,\lambda), I^{(0)}_\epsilon(t,\lambda)) +
(\epsilon|\beta)^2\frac{1}{2(\lambda-u_i)}. 
\een
Therefore the above identities are equivalent to
\ben
\partial_\lambda\log \widetilde{b}_\epsilon(\lambda) =
-\frac{1}{n-2} \operatorname{rk}(\epsilon)^2\lambda^{-1}.
\een
It remain only to recall the definition of
$\widetilde{b}_\epsilon(\lambda)$.

Finally, let us prove that $A'/A''=-1$. Let us compute
\beq\label{c-ratio}
\log (c^i_{\epsilon''}(t,\lambda)/c^i_{\epsilon'}(t,\lambda)) =
\lim_{x\to 0}\ \left(
\int_{
t-\lambda\mathbf{1}}^{t-(u+x)\mathbf{1}} \W_{\epsilon',\epsilon'} + 
\int^{
t-\lambda\mathbf{1}}_{t-(u+x)\mathbf{1}} \W_{\epsilon'',\epsilon''}\right). 
\eeq
Let $C_x$ be a small loop based at $t-(u+x)\mathbf{1}$ that goes once
in a counterclockwise direction around $t-u\mathbf{1}$. Note that if
we add to the RHS of \eqref{c-ratio}  the integral $\lim_{x\to 0}
\int_{C_x}\W_{\epsilon',\epsilon'} $, then the resulting sum of 3
integrals coincides with $\oint_C \W_{\epsilon',\epsilon'} $, where
$C$ is a simple loop  based at
$t-\lambda\mathbf{1}$ that goes around $t-u\mathbf{1}$. We get 
\ben
\log (c^i_{\epsilon''}(t,\lambda)/c^i_{\epsilon'}(t,\lambda)) =
\oint_C \W_{\epsilon',\epsilon'} - \lim_{x\to 0} \int_{C_x}\W_{\epsilon',\epsilon'},
\een
Recalling formula \eqref{phase-form-periods} we get that
\ben
\oint \W_{\epsilon',\epsilon'} =
\lim_{\mu\to \lambda} \Big(
\Omega_{\epsilon'',\epsilon''}(t,\lambda,\mu) -
\Omega_{\epsilon',\epsilon'}(t,\lambda,\mu)\Big)
-2\pi\ii (\epsilon'|\beta)\langle\beta, \epsilon'\rangle .
\een
Note that (the computation is the same as in the proof of
\eqref{phase-period}) 
\ben
\lim_{x\to 0} \int_{C_x}\W_{\epsilon',\epsilon'} =
-(\epsilon'|\beta)^2 \pi\ii=-\pi\ii. 
\een
Recalling also formula \eqref{ope-b} we get
\ben
\frac{b_{\epsilon'}(t,\lambda)}{ b_{\epsilon''}(t,\lambda)} =
\lim_{\mu\to \lambda} \
\exp\Big(
\Omega_{\epsilon',\epsilon'}(t,\lambda,\mu)-
\Omega_{\epsilon'',\epsilon''}(t,\lambda,\mu) +
2\pi\ii\langle \epsilon''-\epsilon', \epsilon^3_1\rangle
\Big).
\een
Combining the above formulas, we get 
$
A'/A'' =- e^{2\pi\ii (
  \langle \epsilon''-\epsilon', \epsilon^3_1\rangle -
  (\epsilon'|\beta)\langle\beta, \epsilon'\rangle )}.
$
We need only to check that the number
\beq\label{ratio-a}
\langle \epsilon''-\epsilon', \epsilon^3_1\rangle -
(\epsilon'|\beta)\langle\beta, \epsilon'\rangle
\eeq
is an integer. Since $(\epsilon'|\beta)=1$ and $\langle
\epsilon',\beta\rangle= 
\frac{1}{2} -\langle \epsilon',\epsilon''\rangle$, we get that
\eqref{ratio-a} coincides with
\beq\label{ratio-aa}
\langle \epsilon''-\epsilon', \epsilon^3_1\rangle -
\langle \epsilon',\epsilon''\rangle -\frac{1}{2}.
\eeq
If $\epsilon''=\pm \epsilon^3_1$, then $\epsilon'\neq \pm
\epsilon^3_1$ and using the explicit formulas for the Euler pairing
(see Section \ref{sec:eup}), we get $\langle
\epsilon',\epsilon^3_1\rangle = \langle
\epsilon',\epsilon''\rangle  = 0$ and $\langle
\epsilon'',\epsilon^3_1\rangle =\pm\tfrac{1}{2}$. Clearly,
\eqref{ratio-aa} is an integer in this case. If $\epsilon'=\pm
\epsilon^3_1$ the argument is identical. Let us assume that both
$\epsilon'$ and $\epsilon''$ are not proportional to
$\epsilon^3_1$. Recalling again the explicit formulas for the Euler
pairing, we get
$
\langle \epsilon',\epsilon^3_1\rangle =
\langle \epsilon'',\epsilon^3_1\rangle  = 0$ and
$
\langle \epsilon',\epsilon''\rangle =\pm\tfrac{1}{2}.$ Clearly, 
\eqref{ratio-aa} is an integer in this case too. 
\qed

\bibliographystyle{amsalpha}

\begin{thebibliography}{FKRW}

\bibitem{CM}
{J. P. Cheng and T. Milanov.}
\textit{The extended D-Toda hierarchy}. Selecta Math. (N.S.) 27 (2021), no. 2, Paper No. 24, 85 pp. 

\bibitem{Du}
{B.~Dubrovin.}
\textit{Geometry of 2D topological field theories}. 
In: ``Integrable systems and quantum groups'' 
(Montecatini Terme, 1993), 120--348, Lecture Notes
in Math., 1620, Springer, Berlin, 1996.

\bibitem{Du2}{B.~Dubrovin.}
\emph{
Painlev\'e transcendents in two dimensional topological field theory.}
arXiv: math/9803107.

\bibitem{G1}
{A.~Givental.}
\textit{Semisimple Frobenius structures at higher genus}. 
Internat. Math. Res. Notices, vol. 23(2001): 1265--1286.

\bibitem{G2}
{A.~Givental.}
\textit{Gromov--Witten invariants and quantization of quadratic Hamiltonians}. 
Mosc. Math. J., vol. 1(2001), 551--568. 

\bibitem{G3}
{A.~Givental.}
\textit{$A_{n-1}$-singularities and $n$-KdV hierarchies.}
Mosc. Math. J., vol. 3, No. 2(2003): 475--505.

\bibitem{GM}
  {A.~Givental and T.~Milanov.}
  \textit{ Simple singularities and integrable hierarchies.}
  The breadth of symplectic and Poisson geometry, 173–-201,
  Progr. Math., 232, Birkhäuser Boston, Boston, MA, 2005.

\bibitem{Ir}{H.~Iritani.}
\emph{An integral structure in quantum cohomology and mirror symmetry
  for toric orbifolds.} Adv. Math., vol. 222, No. 3(2009): 1016–-1079.

\bibitem{Kac}
V.G.~Kac,
\textit{Infinite-dimensional Lie algebras}. 3rd ed.,
Cambridge Univ. Press, Cambridge, 1990. 

\bibitem{Ko}{M. Kontsevich}.
  \textit{Intersection theory on the moduli space of curves and the
    matrix Airy function.}
  Commun. Math. Phys., vol. 147(1992): 1--23.

 \bibitem{KvL}
F.~ ten Kroode and J. van de Leur.
\textit{
Bosonic and fermionic realization of the affine algebra
$\widehat{so}_{2n}$.}
Comm. in Alg., vol. 20, no. 11(1992): 3119--3162.


\bibitem{Mal}{B. Malgrange}.
\emph{Sur les deformations isomonodromiques I: Singularites
  regulieres. }(1982)


\bibitem{Mi}{T. Milanov}.
  \textit{The Period map for quantum cohomology of $\mathbb{P}^2$.}
  Adv. in Math., vol. 351(2019): 804–-869.
  

\bibitem{Mi2}{T. Milanov}.
\textit{The phase factors in singularity theory.} In: ``Primitive forms and related subjects-Kavli IPMU 2014", Adv. Studies in Pure Math, Vol. 83, Math. Soc of Japan, Tokyo, 2019, pp. 295--326.

\bibitem{Mi3}{T. Milanov}. 
\textit{Lectures on Painleve property for semi-simple Frobenius
  manifolds.}
Preprint arXiv:1702.04323

\bibitem{Mi4}{T. Milanov}.
  \textit{The total ancestor potential in singularity theory.}
  In: ``B-model Gromov--Witten theory'', 
  Trends in Mathematics,
  E. Clauder and Y. Ruan (eds.),
  Springer Nature Switzerland AG, 2018, pp. 539--571.

\bibitem{MST}{T. Milanov, Y. Shen, and H.-H. Tseng.}
\textit{Gromov--Witten theory of Fano orbifold curves, Gamma integral
  structures and ADE-Toda hierarchies.}
Geometry \& Topology, vol. 20(2016): 2135--2218.

\bibitem{MT}{T. Milanov and H.-H. Tseng}.
  \textit{The spaces of Laurent polynomials, Gromov--Witten theory of
    $\PP^1$-orbifolds, and integrable hierarchies.}
  J. Reine Angew. Math., vol. 622(2008): 189--235.


\bibitem{Te}
C.~Teleman.
\emph{The structure of 2D semi-simple field theories.} 
Invent. Math., vol. 188, no. 3(2012): 525--588.

\bibitem{Wi}{E. Witten}.
  \textit{
    Two-dimensional gravity and intersection theory on moduli space.}
  Surv. Diff. Geom., vol. 1(1991): 243--310.

\end{thebibliography}

\end{document}